\documentclass[11pt]{amsart}
\usepackage[utf8]{inputenc}
\usepackage[T1]{fontenc}
\usepackage{lmodern}
\usepackage{amssymb,amsmath,amsthm,mathtools}
\usepackage{booktabs,array}
\usepackage[margin=1.15in]{geometry}
\usepackage{enumitem}
\usepackage[dvipsnames]{xcolor}
\usepackage[colorlinks=true,linkcolor=blue!50!black,citecolor=green!35!black,urlcolor=blue!50!black]{hyperref}
\hypersetup{pdftitle={The generic splitting type of the logarithmic bundle of an arrangement of at most sixteen lines is combinatorial},pdfauthor={Cristian Anghel, Rares Stan}}

\newtheorem{theorem}{Theorem}[section]
\newtheorem{proposition}[theorem]{Proposition}
\newtheorem{lemma}[theorem]{Lemma}
\newtheorem{corollary}[theorem]{Corollary}
\newtheorem*{theoremA}{Theorem A}
\newtheorem*{theoremB}{Theorem B}
\newtheorem*{theoremC}{Theorem C}
\theoremstyle{definition}
\newtheorem{definition}[theorem]{Definition}
\newtheorem{conjecture}[theorem]{Conjecture}
\newtheorem{question}[theorem]{Question}
\theoremstyle{remark}
\newtheorem{remark}[theorem]{Remark}
\newtheorem{example}[theorem]{Example}

\newcommand{\cA}{\mathcal{A}}
\newcommand{\cB}{\mathcal{B}}
\newcommand{\cC}{\mathcal{C}}
\newcommand{\cF}{\mathcal{F}}
\newcommand{\cI}{\mathcal{I}}
\newcommand{\cO}{\mathcal{O}}
\newcommand{\cR}{\mathcal{R}}
\newcommand{\CC}{\mathbb{C}}
\newcommand{\QQ}{\mathbb{Q}}
\newcommand{\PP}{\mathbb{P}}
\newcommand{\mm}{\mathfrak{m}}
\newcommand{\mdr}{\operatorname{mdr}}
\newcommand{\Der}{\operatorname{Der}}
\newcommand{\len}{\operatorname{len}}
\newcommand{\ord}{\operatorname{ord}}
\newcommand{\mult}{\operatorname{mult}}
\newcommand{\Sing}{\operatorname{Sing}}
\newcommand{\inn}{\operatorname{in}}
\newcommand{\TA}{\mathcal{T}\langle\cA\rangle}
\newcommand{\GST}{\textup{(GST)}}

\title[Generic splitting type for at most sixteen lines]{The generic splitting type of the logarithmic bundle\\ of an arrangement of at most sixteen lines\\ is combinatorial}

\author{Cristian Anghel}
\address{Simion Stoilow Institute of Mathematics of the Romanian Academy, P.O.\ Box 1-764, RO-014700 Bucharest, Romania}
\email{cristian.anghel@imar.ro}

\author{Rareș Stan}
\address{Simion Stoilow Institute of Mathematics of the Romanian Academy, P.O.\ Box 1-764, RO-014700 Bucharest, Romania}
\email{rstan@imar.ro}


\begin{document}

\begin{abstract}
Let $\cA$ be an arrangement of $d$ lines in $\PP^2_\CC$ and let $\TA$ be its bundle of logarithmic vector fields. On a general line, $\TA$ splits as $\cO(-k)\oplus\cO(k+1-d)$ with $k=\min\{\mdr(\cA),\lfloor\frac{d-1}{2}\rfloor\}$, where $\mdr$ is the minimal degree of a Jacobian syzygy. We prove that $k$ is determined by the intersection lattice whenever $d\le16$, answering in this range a question of Cook, Harbourne, Migliore and Nagel. Since freeness can be read off from $k$ and $c_2(\TA)$, this strengthens Terao's freeness conjecture in this range; in particular we obtain a proof of Terao's conjecture for $14$ lines that uses it only for $13$ lines. For any number $d$ of lines, it follows that $\mdr$ is determined by the lattice whenever $\mdr\le6$ and $2\mdr\le d-3$, and that Terao's conjecture holds for every free arrangement whose smaller exponent is at most $7$.

The proof reduces the question to \emph{critical} arrangements: $2s+3$ lines, each containing at most $s+1$ intersection points, with $\mdr=s$. For these, the inflection polynomial of the minimal derivation factors as $Q_\cA R$ with $R\in H^0(\cI_Z^2(s-3))$, where $Z$ is the zero scheme of the derivation, of length $\nu=c_2(\TA)-s(s+2)$. For $s\le5$ this forces $\nu\le1$, and Chern class arguments then reduce everything to Terao's conjecture for $13$ lines. For $s=6$ we show $\nu\le3$, using a Riccati equation along multiple components of $R$, Ziegler's multirestriction, and local computations verified by computer algebra; an exhaustive enumeration of the $26\,265\,318$ candidate lattices with $\nu\in\{2,3\}$ that satisfy Bojanowski's inequality completes the proof.
\end{abstract}

\maketitle
\tableofcontents

%======================================================================
\section{Introduction}
%======================================================================

Let $\cA=\{H_1,\dots,H_d\}$ be an arrangement of $d$ distinct lines in $\PP^2=\PP^2_\CC$, $H_i=\{\alpha_i=0\}$, $Q_\cA=\prod\alpha_i$, and let $S=\CC[x,y,z]$. The module of logarithmic derivations $D(\cA)=\{\theta\in\Der S:\theta(\alpha_i)\in\alpha_iS\}$ splits as $SE\oplus D_0(\cA)$, where $E$ is the Euler derivation and $D_0(\cA)$ is the module of Jacobian syzygies of $Q_\cA$. Its sheafification $\TA$ is a rank-$2$ vector bundle on $\PP^2$, the \emph{logarithmic bundle}, with
\[
c_1(\TA)=1-d,\qquad c_2(\TA)=\tbinom{d-1}{2}-X(\cA),\qquad X(\cA)=\sum_p\tbinom{m_p-1}{2},
\]
the sum running over the intersection points $p$, of multiplicity $m_p$. The arrangement is \emph{free} if $D_0(\cA)$ is a free $S$-module, i.e.\ if $\TA$ splits. Terao's conjecture \cite{OT92} asserts that freeness depends only on the intersection lattice $L(\cA)$. It is known for at most $14$ lines (Faenzi--Vallès \cite{FV14} up to $12$, Dimca--Ibadula--Măcinic \cite{DIM19} for $13$, Barakat et al.\ \cite{BBJKL21,BK23} for $14$, for rank-$3$ arrangements in any characteristic) and for $15$ and $16$ lines \cite{AS}.

Let $\mdr(\cA)$ be the smallest degree of a nonzero element of $D_0(\cA)$. On a general line $\ell$ one has (Proposition~\ref{prop:gst})
\[
\TA|_\ell\cong\cO_\ell(-k)\oplus\cO_\ell(k+1-d),\qquad k=k(\cA):=\min\Bigl\{\mdr(\cA),\Bigl\lfloor\tfrac{d-1}2\Bigr\rfloor\Bigr\}.
\]
We call $k(\cA)$ the \emph{generic splitting index}. Since $\cA$ is free if and only if $c_2(\TA)=k(d-1-k)$ (Proposition~\ref{prop:freecrit}), the following conjecture implies Terao's conjecture for line arrangements.

\begin{conjecture}[GST]\label{conj:gst}
If $L(\cA)\cong L(\cB)$, then $k(\cA)=k(\cB)$; i.e.\ the generic splitting type of the logarithmic bundle is combinatorial.
\end{conjecture}

This is the question asked by Cook, Harbourne, Migliore and Nagel \cite[Question~7.12]{CHMN18}, stated as a conjecture by Dimca \cite[Conj.~3.9]{Dim24}; see also \cite{DP26}. Dimca also conjectured that the maximal dimension of a graded piece of the Jacobian module $I/J$, where $J$ is the Jacobian ideal of $Q_\cA$ and $I$ its saturation, is determined by the lattice \cite[Conj.~3.7]{Dim24}. Dimca and Pokora recently studied this statement as a generalized Terao conjecture \cite[Conj.~1.1]{DPgen}. They observed that it is equivalent to a positive answer to \cite[Question~7.12]{CHMN18}, i.e.\ to \GST, and proved it for at most $12$ lines and for $13$ lines whose maximal multiplicity is not $5$ \cite[Thm.~1.3]{DPgen}. Theorems~A and~B below extend this to all arrangements of at most $16$ lines. Their deletion lemma \cite[Lemma~3.4]{DPgen} is of the same kind as Theorem~\ref{thm:structure}(3), on which our reduction formula is based.

The integer $\mdr$ itself is not combinatorial: for Ziegler's pair of $9$ lines with $6$ triple points \cite{Zie89a}, $\mdr$ is $5$ or $6$ according to whether the triple points lie on a conic. But in both cases $k=4$, as \GST\ predicts. In the dual plane, $k(\cA)$ governs the curves through the dual configuration of points that have a point of high multiplicity at a general point, by the work of Cook, Harbourne, Migliore and Nagel \cite{CHMN18}. For the critical arrangements that carry the whole difficulty, \GST\ asks whether the existence of an \emph{unexpected curve} of a given degree depends only on the matroid (Remark~\ref{rem:dual}).

\begin{theoremA}
\GST\ holds for all arrangements of at most $14$ lines. The proof uses Terao's conjecture only for free arrangements of $13$ lines with exponents $(6,6)$ \cite{DIM19}, and no computer assistance; as far as we can see, the proof in \cite{DIM19} does not use computers either. In particular, Terao's conjecture for $14$ lines follows from Terao's conjecture for $13$ lines.
\end{theoremA}

\begin{theoremB}
\GST\ holds for all arrangements of at most $16$ lines. In addition to Theorem~A, the proof uses Terao's conjecture only for free arrangements of $15$ lines with exponents $(7,7)$ \cite{AS}, computer-verified local computations, and an exhaustive enumeration of rank-$3$ matroids.
\end{theoremB}

Since the reduction to critical arrangements works for any number of lines, Theorems~A and~B also have consequences for arrangements of more than $16$ lines.

\begin{theoremC}
Let $\cA$ and $\cB$ be arrangements of $d$ lines with $L(\cA)\cong L(\cB)$, where $d$ is arbitrary.
\begin{enumerate}[label=\textup{(\alph*)}]
\item If $\min\{k(\cA),k(\cB)\}\le6$, then $k(\cA)=k(\cB)$. In particular, if $\mdr(\cA)\le6$ and $2\mdr(\cA)\le d-3$, then $\mdr(\cB)=\mdr(\cA)$.
\item If $\cA$ is free with exponents $(d_1,d_2)$, $d_1\le d_2$, and $d_1\le7$, then $\cB$ is free.
\end{enumerate}
Using only Theorem A, the same holds with $6$ and $7$ replaced by $5$ and $6$; this version does not use computer assistance.
\end{theoremC}

Part (b) extends the case $d_1\le5$ of Terao's conjecture, due to Abe \cite[Cor.~5.5]{Abe14}, to $d_1\le7$. Since $2\mdr(\cA)\le d-3$ means $\mdr(\cA)<(d-2)/2$, part (a) proves, for $\mdr\le6$, the conjecture \cite[Conj.~1.2]{DP26} that $\mdr$ is determined by the lattice when $\mdr<(d-2)/2$. By \cite{DP26}, it is equivalent to \cite[Conj.~3.7]{Dim24}, the conjecture recalled above. It is also equivalent to \GST\ directly, since $\mdr<(d-2)/2$ holds if and only if $\mdr\le K-1$, where $K=\lfloor\frac{d-1}2\rfloor$, i.e.\ if and only if $k=\mdr<K$; so Theorems A and B prove it for $d\le16$. It is a weaker form of Dimca's conjecture \cite[Conj.~3.5]{Dim24}, with $\mdr<d/2$, which Dimca and Pokora disprove in \cite{DP26} with two arrangements of $9$ lines with the same lattice and $\mdr$ equal to $4$ and $5$. The condition $2\mdr(\cA)\le d-3$ in (a) therefore cannot be relaxed to $\mdr(\cA)<d/2$, and it cannot be dropped either, as Ziegler's pair shows; in both pairs the two values of $k$ agree, as \GST\ predicts. By Remark~\ref{rem:dual}, part (a) says the following for two configurations of points in $\PP^2$ with isomorphic matroids, of any size: if for one of them the least $j$ such that some curve of degree $j+1$ through the points has multiplicity at least $j$ at a general point is at most $6$, then it is the same for the other. Finally, Terao's conjecture for $17$ and $18$ lines becomes a part of the next case $(C_7)$ of our method (Proposition~\ref{prop:1718}).

\subsection*{Strategy}
\emph{Reduction (\S\ref{sec:reduction}).} Deletion--restriction arguments show that $k(\cA)$ is computed from the lattice together with the predicates $[\mdr\le s]$ on sub-arrangements of a very special kind, the \emph{$s$-candidates}: $2s+3$ lines, each containing at most $s+1$ intersection points (Theorem~\ref{thm:reduction}). Hence \GST\ for $d\le2S+4$ follows from the following condition for all $s\le S$:
\begin{quote}
$(C_s)$: for every $s$-candidate lattice all of whose realizations have $\mdr\ge s$, either all realizations or none have $\mdr=s$.
\end{quote}
A realization with $\mdr=s$ is called \emph{critical}.

\emph{The canonical foliation (\S\ref{sec:foliation}).} The minimal derivation $\theta$ of a critical arrangement defines a saturated foliation $\cF$ of degree $s$ for which all lines of $\cA$ are invariant. Its inflection polynomial factors as $F_\theta=Q_\cA R$ with $\deg R=s-3$. The \emph{residual lemma} (Theorem~\ref{thm:residual}) says that $R$ lies in the square of the ideal of the zero scheme $Z$ of $\theta$, a scheme of length
\[
\nu=c_2(\TA)-s(s+2),
\]
the \emph{defect}. Thus $\nu=0$ if and only if the critical arrangement is free. (This defect should not be confused with the invariant denoted $\nu(\cA)$ in \cite{Dim24,DPgen}, the maximal dimension of a graded piece of the Jacobian module.)

\emph{Small defect (\S\ref{sec:nu1}).} Consider a lattice all of whose realizations have $\mdr\ge s$. If a critical realization has $\nu=0$, Bogomolov's inequality shows that every realization is critical. If $\nu=1$, a Chern class argument (Lemma~\ref{lem:c1c2zero}) shows that any other realization with $\mdr\ge s+1$ is free with exponents $(s+1,s+1)$. Hence Terao's conjecture for these free arrangements of $2s+3$ lines implies $(C_s)$ for such lattices (Proposition~\ref{prop:nu1}). For $s\le5$ the residual lemma forces $\nu\le1$ (Theorem~\ref{thm:nu13}); this proves Theorem A.

\emph{Fifteen lines (\S\ref{sec:fifteen}).} For $s=6$, $R$ is a cubic, and $\nu\le1$ fails: the $15$-line subsets of the reflection arrangement $\cA(G(5,1,3))$ include critical arrangements with $\nu=2$. We prove $\nu\le3$ (Theorem~\ref{thm:nu15}). If $R$ is not the cube of a linear form, this follows from local lemmas on the initial forms of $R$. If $R=L^3$, the line $\ell=\{L=0\}$ carries all of $Z$ and three mechanisms take over:
\begin{itemize}[leftmargin=*]
\item a \emph{Riccati equation} along $\ell$ when $\ell$ is invariant: the leaves near $\ell$ osculate a pencil of lines centred at a point $p_0\in\ell$, and the tangential restriction of $\theta$ divides a product of linear forms;
\item \emph{Ziegler's multirestriction} when $\ell$ is invariant but not in $\cA$;
\item a \emph{degree budget} for the normal component of $\theta$ along $\ell$ when $\ell$ is not invariant.
\end{itemize}
The local input consists of finitely many statements about the jets of $\theta$ at points of $\ell$. They were verified by exact computer algebra, up to the orders allowed by the budgets (Appendix~\ref{app:local}).

\emph{Enumeration (\S\ref{sec:enum}).} By $\nu=43-X(\cA)$, it remains to treat the $6$-candidate lattices with $X\in\{40,41\}$; Bojanowski's inequality, valid for all complex line arrangements of this kind, cuts them down further. We enumerated them with the programs of \cite{AS}: there are $26\,265\,318$ isomorphism classes. Exactly six are realizable over $\CC$. Two of them are critical, and each has a single Galois orbit of realizations, so $(C_6)$ holds.

\emph{Any number of lines (\S\ref{sec:anyd}).} The reduction formula holds for every $d$. Hence $(C_s)$ for $s\le6$ determines $[\mdr\le s]$ from the lattice for every $s\le6$ with $s<K$, whatever $d$ is. For a free arrangement with $d_1\le7$, this gives $\mdr\ge d_1$ for every arrangement with the same lattice, and $c_2=d_1d_2$ then forces freeness. This gives Theorem C, and shows that Terao's conjecture for $17$ lines is equivalent to the stratum $\nu=1$ of $(C_7)$, while Terao's conjecture for $18$ lines follows from the strata $\nu\in\{2,3\}$.

\subsection*{Computer assistance}
Two kinds of computation enter the proof of Theorem B, both in exact arithmetic:
\begin{itemize}[leftmargin=*]
\item local computations with truncated power series over $\QQ[\text{parameters}]$, used in \S\ref{sec:fifteen};
\item the enumeration and realizability screening of \S\ref{sec:enum}, which reuses the programs validated in \cite{AS}.
\end{itemize}
A rank computation for multiarrangements on $\PP^1$ gives an independent check of Lemma~\ref{lem:multiarr}, which follows from a theorem of Abe.
Theorem A does not depend on any of them. Theorem B also depends on \cite{AS} in two ways: through its main result, Terao's conjecture for $15$ lines, which is itself computer-assisted and is used only for free arrangements with exponents $(7,7)$, and through its generation and screening programs. Theorem C and Proposition~\ref{prop:1718} depend on Theorem B, hence on these computations and on \cite{AS}; the version of Theorem C with $5$ and $6$ depends only on Theorem A. At the time of writing \cite{AS} is a preprint. The programs, the screening logs, the realizable classes and SHA-256 checksums of the complete enumeration output (about $2.9$\,GB, available on request and deterministically regenerable with the programs) are provided as ancillary files.

\subsection*{Conventions}
$\cA$ always denotes an arrangement of $d\ge3$ lines, not a pencil. For $H\in\cA$ we write $\cA^H$ for the set of intersection points on $H$ and $r_H=|\cA^H|$; thus $\sum_{p\in\cA^H}(m_p-1)=d-1$. We put $X=\sum_p\binom{m_p-1}2$, $W=\sum_p\binom{m_p-2}2$, $\tau=\sum_p(m_p-1)^2$ and $K=\lfloor\frac{d-1}2\rfloor$. Degrees of derivations are degrees of their coefficients. The results of \S\ref{sec:prelim} also hold for pencils, and \S\ref{sec:reduction} applies them to sub-arrangements that may be pencils. In the local computations of \S\S\ref{sec:foliation}, \ref{sec:nu1}, \ref{sec:fifteen} and Appendix~\ref{app:local}, the letters $a,b,c,f,g,A,B$ denote local coefficients and are redefined in each computation, and $X$ without argument denotes the affine vector field of a derivation; there the combinatorial invariant is always written $X(\cA)$.

%======================================================================
\section{Preliminaries}\label{sec:prelim}
%======================================================================

\subsection{The logarithmic bundle}
$\TA$ is reflexive of rank $2$ on a smooth surface, hence locally free, and $D_0(\cA)=H^0_*(\TA)$. Its Chern classes are $c_1=1-d$ and $c_2=(d-1)^2-\tau=b_2-(d-1)$ with $b_2=\sum_p(m_p-1)$. Since $\sum_p\binom{m_p}2=\binom d2$, this gives $c_2=\binom{d-1}2-X$. Moreover $\TA(1)\cong T_{\PP^2}(-\log\cA)$, the sheaf of vector fields tangent to $\cA$.

\begin{proposition}[generic splitting; cf.\ \cite{AD18}]\label{prop:gst}
For a general line $\ell$, $\TA|_\ell\cong\cO_\ell(-k)\oplus\cO_\ell(k+1-d)$ with $k=k(\cA)$.
\end{proposition}
\begin{proof}
Let $r=\mdr(\cA)$ and $0\neq\theta\in D_0(\cA)_r$.

\emph{Case $2r<d-1$.} The section $\theta$ of $\TA(r)$ has no divisorial zeros: otherwise dividing by the equation of the divisor would give an element of $D_0(\cA)$ of smaller degree. This gives an exact sequence $0\to\cO(-r)\to\TA\to\cI_Z(r+1-d)\to0$ with $Z$ finite. On a line avoiding $Z$ the induced extension splits, since $H^1(\cO_\ell(d-1-2r))=0$.

\emph{Case $2r\ge d-1$.} Then $\TA$ is semistable, since $H^0(\TA(t))=D_0(\cA)_t$ vanishes whenever $c_1(\TA(t))=1-d+2t<0$; by Grauert--Mülich \cite{OSS} its generic splitting type is as balanced as possible.
\end{proof}

The following is well known; we include the short proof.

\begin{proposition}\label{prop:freecrit}
$\cA$ is free if and only if $c_2(\TA)=k(d-1-k)$. If $2\mdr(\cA)<d-1$, then $c_2(\TA)=\mdr(\cA)(d-1-\mdr(\cA))+\len Z$, with $Z$ as in the proof of Proposition~\ref{prop:gst}.
\end{proposition}
\begin{proof}
The second statement follows from the exact sequence. If $2r<d-1$ and $c_2=r(d-1-r)$, then $Z=\varnothing$ and the sequence splits. If $2r\ge d-1$, Bogomolov's inequality for the semistable bundle gives $(d-1)^2\le 4c_2$. For $d$ even this excludes $c_2=k(d-1-k)$. For $d$ odd the normalized bundle is semistable with $c_1=c_2=0$, hence trivial. The converse is clear.
\end{proof}

\begin{corollary}
\GST\ for arrangements of $d$ lines implies Terao's conjecture for $d$ lines.
\end{corollary}

\begin{remark}[dual form]\label{rem:dual}
Let $\ell$ be a general line. In the dual plane, the lines of $\cA$ become a set $\cA^\vee$ of $d$ points and $\ell$ becomes a general point $P$. By \cite[Lemma~3.3]{CHMN18}, the forms of degree $j+1$ vanishing on $\cA^\vee$ and with multiplicity $\ge j$ at $P$ form a space of dimension $\max(0,j-k+1)+\max(0,j-d+2+k)$, where $k=k(\cA)$. So $k(\cA)$ is the least $j$ for which such a curve exists.

Now let $\cA$ realize an $s$-candidate, so $d=2s+3$. Every multiplicity is $\le s+1$ (Lemma~\ref{lem:cand}(a)), so no $s+3$ points of $\cA^\vee$ are collinear. Since $|\cA^\vee|=2(s+1)+1$, the classical lemma on at most $2j+1$ points (see e.g.\ \cite{EGH96}), applied in degree $j=s+1$, then shows that $\cA^\vee$ imposes independent conditions on forms of degree $s+1$. The expected dimension is then $\binom{s+3}2-\binom{s+1}2-(2s+3)=0$. Hence $\mdr(\cA)\le s$ if and only if $\cA^\vee$ has an unexpected curve of degree $s+1$ in the sense of \cite{CHMN18}. Condition $(C_s)$ of Theorem~\ref{thm:induction} therefore says that, for these lattices, the existence of such a curve depends only on the matroid.

We do not claim that \GST\ is equivalent to the combinatorial invariance of unexpected curves in general. That notion also involves the Hilbert function of $\cA^\vee$, which is not determined by the matroid; consider, for example, six points on a conic.
\end{remark}

\subsection{Restriction and deletion}
For $\theta\in D(\cA)_e$ and a line $H$ we write $\theta|_H\equiv0$ if $\theta(p)\in\CC p$ for every $p\in H$, i.e.\ if $\alpha_H$ divides the coefficients of $\theta\wedge E$. If $H$ is $\theta$-invariant and $\theta|_H\not\equiv0$, then $\theta$ induces a nonzero vector field on $H\cong\PP^1$ with values in $\cO_H(e-1)$, i.e.\ a section
\[
t_H(\theta)\in H^0(\cO_H(e+1)),
\]
the \emph{tangential restriction}. For an arbitrary line $\ell$, the \emph{normal restriction} is $n_\ell(\theta)=\theta(\alpha_\ell)|_\ell\in H^0(\cO_\ell(e))$. It does not change under $\theta\mapsto\theta+hE$, it vanishes identically if and only if $\ell$ is $\theta$-invariant, and otherwise its zeros are the points of $\ell$ at which the foliation is singular or tangent to $\ell$.

\begin{lemma}\label{lem:restr}
Let $H$ be a line, and let $\cB$ be an arrangement not containing $H$.
\begin{enumerate}[label=\textup{(\alph*)}]
\item If $H\in\cA$, $\theta\in D(\cA)_e$ and $\theta|_H\not\equiv0$, then $t_H(\theta)$ vanishes at every point of $\cA^H$; in particular $r_H\le e+1$.
\item If $\psi\in D(\cB)_{e}$, then $\ord_p n_H(\psi)\ge m_p(\cB)-1$ at every point $p\in H$ lying on $m_p(\cB)\ge1$ lines of $\cB$. Hence, if $n_H(\psi)\ne0$, then $\sum_{p\in H}(m_p(\cB)-1)=|\cB|-|\cB\cap H|\le e$, where $\cB\cap H$ is the set of points of $H$ lying on lines of $\cB$.
\end{enumerate}
\end{lemma}
\begin{proof}
(a) At $p\in\cA^H$, $\theta(p)$ is tangent to two distinct lines, so it is proportional to $p$.

(b) Take affine coordinates with $p=0$, $H=\{y=0\}$, the lines of $\cB$ through $p$ being $y=c_jx$, $1\le j\le m$, with distinct $c_j$. Write $\psi=a\partial_x+b\partial_y$ and let $a_i,b_i$ be the homogeneous parts of degree $i$. Tangency gives $\varphi_i(c):=b_i(1,c)-c\,a_i(1,c)=0$ for $c=c_j$. Since $\deg\varphi_i\le i+1$, we get $\varphi_i\equiv0$ for $i\le m-2$, hence $b_i(1,0)=0$. Thus $n_H(\psi)=b(x,0)$ has order $\ge m-1$ at $p$.
\end{proof}

\begin{lemma}\label{lem:factor}
If $\theta\in D(\cA)_e$ and $\theta|_H\equiv0$ for every $H$ in a subset $G\subseteq\cA$, then $\theta=g\,\theta_0+hE$ with $g=\prod_{H\in G}\alpha_H$, $\theta_0\in D(\cA\setminus G)_{e-|G|}$ and $h\in S$. If $\theta\notin SE$, then $\theta_0\notin SE$.
\end{lemma}
\begin{proof}
$\theta\wedge E$ is divisible by $g$. Since $S/(g)$ has depth $2$, the Koszul cohomology $H^1(x,y,z;S/(g))$ vanishes, i.e.\ the complex $S/(g)\to(S/(g))^3\to\bigwedge^2(S/(g))^3$ given by $\wedge E$ is exact in the middle. Hence $\theta\equiv hE$ modulo $g$. Tangency of $\theta_0$ to the lines outside $G$ is clear.
\end{proof}

\begin{corollary}[gap lemma]\label{cor:gap}
Let $r=\mdr(\cA)$ and $0\ne\theta\in D_0(\cA)_r$. For every $H\in\cA$, either $\theta|_H\not\equiv0$ and $r_H\le r+1$, or $\theta|_H\equiv0$ and $r_H\ge d-r$.
\end{corollary}
\begin{proof}
If $\theta|_H\equiv0$, Lemma~\ref{lem:factor} gives $\theta=\alpha_H\theta'+hE$ with $\theta'\in D(\cA\setminus H)_{r-1}\setminus SE$. Then $\theta'(\alpha_H)\notin(\alpha_H)$, since otherwise $\theta'$ would give an element of $D_0(\cA)$ of degree $r-1$. Now apply Lemma~\ref{lem:restr}(b) to $\cB=\cA\setminus H$ and $e=r-1$: we get $d-1-r_H\le r-1$.
\end{proof}

\begin{lemma}[pencil bound]\label{lem:pencil}
If $\cA$ has a point of multiplicity $m$, then $\mdr(\cA)\le d-m$.
\end{lemma}
\begin{proof}
If $p=(0{:}0{:}1)$, then $\bigl(\prod_{H\not\ni p}\alpha_H\bigr)\partial_z\in D(\cA)\setminus SE$.
\end{proof}

\begin{proposition}[addition--deletion; cf.\ \cite{ADS21}]\label{prop:ad}
Let $H\in\cA$, $\cA'=\cA\setminus H$, $r=\mdr(\cA)$ and $r'=\mdr(\cA')$.
\begin{enumerate}[label=\textup{(\roman*)}]
\item $r'\le r\le r'+1$.
\item If $d-1-r_H>r'$, then $r=r'$.
\item If $r=r'$, then $r_H\le r+1$.
\end{enumerate}
\end{proposition}
\begin{proof}
(i) Use $D_0(\cA)\subseteq D(\cA')$ and $\alpha_HD_0(\cA')\subseteq D(\cA)$. (ii) A nonzero $\theta'\in D_0(\cA')_{r'}$ is tangent to $H$ by Lemma~\ref{lem:restr}(b), so $r\le r'$. (iii) If $\theta\in D_0(\cA)_r$ had $\theta|_H\equiv0$, Lemma~\ref{lem:factor} would give $r'\le r-1$; so $\theta|_H\not\equiv0$, and Lemma~\ref{lem:restr}(a) applies.
\end{proof}

%======================================================================
\section{Reduction to critical candidates}\label{sec:reduction}
%======================================================================

In the \emph{unstable range} $2\mdr(\cA)\le d-3$ we have $k=\mdr(\cA)<K$.

\begin{theorem}[structure in the unstable range]\label{thm:structure}
Let $r=\mdr(\cA)$ with $2r\le d-3$. Put $G=\{H\in\cA:r_H\ge r+2\}$, $\cA_0=\cA\setminus G$ and $r_0=r-|G|$. Then:
\begin{enumerate}[label=\textup{(\arabic*)}]
\item $G$ is the set of lines on which the minimal derivation vanishes, and $\theta=\bigl(\prod_G\alpha_H\bigr)\theta_0+hE$;
\item $\mdr(\cA_0)=r_0\ge0$, $|\cA_0|\ge2r_0+3$, and every line of $\cA_0$ contains at most $r_0+1$ points of $\cA_0$;
\item every $\cB\subseteq\cA_0$ with $|\cB|\ge2r_0+2$ has $\mdr(\cB)=r_0$.
\end{enumerate}
For $r_0=0$, $\cA_0$ is a pencil, and \textup{(3)} is trivial.
\end{theorem}
\begin{proof}
(1) Since $r+1<d-r$, this follows from Corollary~\ref{cor:gap} and Lemma~\ref{lem:factor}.

(2) $\theta_0$ gives $\mdr(\cA_0)\le r_0$, and multiplication by $\prod_G\alpha_H$ gives $r\le\mdr(\cA_0)+|G|$. For $H\in\cA_0$ we have $\theta_0|_H\not\equiv0$, so Lemma~\ref{lem:restr}(a) bounds the number of points. Finally $|\cA_0|=d-|G|\ge2r+3-|G|\ge2r_0+3$.

(3) Delete lines one at a time. For $\cC\subseteq\cA_0$ with $c=|\cC|\ge2r_0+3$ and $H\in\cC$, $H$ contains at most $r_0+1$ points of $\cC$. Moreover $\theta_0\in D(\cA_0)\subseteq D(\cC\setminus H)$ and $\theta_0\notin SE$, so $\mdr(\cC\setminus H)\le r_0$. Hence $c-1-|\cC^H|\ge r_0+1>\mdr(\cC\setminus H)$, and Proposition~\ref{prop:ad}(ii) gives $\mdr(\cC)=\mdr(\cC\setminus H)$. Starting from $\cA_0$, this shows $\mdr=r_0$ for all subsets of size $\ge2r_0+2$ obtained by successive deletions, i.e.\ for all of them.
\end{proof}

A deletion principle of the same kind, for pairs of arrangements with isomorphic lattices, is \cite[Lemma~3.4]{DPgen}.

\begin{definition}
An \emph{$s$-candidate} is a simple rank-$3$ matroid on $2s+3$ elements, not a pencil, in which every element lies on at most $s+1$ rank-$2$ flats. A realization $\cA$ of an $s$-candidate is \emph{critical} if $\mdr(\cA)=s$.
\end{definition}

For a simple rank-$3$ matroid $\Lambda$ we write $\cR(\Lambda)$ for the set of arrangements $\cA$ with $L(\cA)\cong\Lambda$, up to the action of $PGL_3$; $\mdr$ and $k$ are constant on $PGL_3$-orbits.

For $s<K$ consider the following condition $\mathrm{P}(s)$, whose first clauses depend only on $L(\cA)$:
\begin{quote}
$G_s=\{H:r_H\ge s+2\}$, $\cA_0=\cA\setminus G_s$, $s_0=s-|G_s|\ge0$; every line of $\cA_0$ contains at most $s_0+1$ points of $\cA_0$; $|\cA_0|\ge2s_0+3$; and $\mdr(\cB)\le s_0$ for one, equivalently every, $(2s_0+3)$-subset $\cB\subseteq\cA_0$.
\end{quote}
Every such $\cB$ is an $s_0$-candidate or a pencil. A pencil has $\mdr=0$, so in that case the last clause holds trivially; this always happens if $s_0=0$, since then every line of $\cA_0$ contains at most one point of $\cA_0$.

\begin{theorem}[reduction formula]\label{thm:reduction}
$k(\cA)=\min\bigl(K,\min\{s<K:\mathrm{P}(s)\}\bigr)$.
\end{theorem}
\begin{proof}
Suppose $\mathrm{P}(s)$ holds. Starting from $\cB$, add the lines of $\cA_0$ one at a time. If $\cC=\cB'\cup\{H\}$ with $|\cB'|\ge2s_0+3$, then $|\cC|-1-|\cC^H|\ge s_0+2>\mdr(\cB')$, so Proposition~\ref{prop:ad}(ii) keeps $\mdr\le s_0$. Hence $\mdr(\cA)\le\mdr(\cA_0)+|G_s|\le s$.

Conversely, if $r=\mdr(\cA)<K$, then $\mathrm{P}(r)$ holds by Theorem~\ref{thm:structure}, and $\mathrm{P}(s)$ fails for $s<r$. If $r\ge K$, then $\mathrm{P}(s)$ fails for all $s<K$. ``One'' equals ``every'' in $\mathrm{P}(s)$ because both are equivalent to $\mdr(\cA_0)\le s_0$: one implication is the first paragraph, and the other holds because $D(\cA_0)\subseteq D(\cB)$.
\end{proof}

\begin{lemma}\label{lem:cand}
Let $\cA$ realize an $s$-candidate.
\begin{enumerate}[label=\textup{(\alph*)}]
\item Every point has multiplicity $\le s+1$.
\item If $\cA$ is critical, then $c_2(\TA)\ge s(s+2)$, with equality if and only if $\cA$ is free; the difference
\[
\nu(\cA):=c_2(\TA)-s(s+2)=\tbinom{2s+2}{2}-s(s+2)-X(\cA)
\]
is a lattice invariant, the \emph{defect}.
\end{enumerate}
\end{lemma}
\begin{proof}
(a) A line not through a point $p$ meets the $m_p$ lines through $p$ in distinct points. (b) Proposition~\ref{prop:freecrit} with $2s<d-1$.
\end{proof}

\begin{theorem}[induction]\label{thm:induction}
Let $S\ge0$ and assume $(C_s)$ for all $s\le S$:
\begin{quote}
$(C_s)$: for every $s$-candidate $\Lambda$ all of whose realizations have $\mdr\ge s$, the set of critical realizations of $\Lambda$ is empty or all of $\cR(\Lambda)$.
\end{quote}
Then $[\mdr\le s]$ is constant on $\cR(\Lambda)$ for every $s$-candidate $\Lambda$ with $s\le S$, and \GST\ holds for all arrangements of at most $2S+4$ lines.
\end{theorem}
\begin{proof}
Induction on $s$. For a realization $\cB$ of an $s$-candidate, $K_\cB=s+1$. By Theorem~\ref{thm:reduction}, $[\mdr(\cB)\le s-1]$ is decided by $L(\cB)$ and by predicates $[\mdr\le s_0]$ on $s_0$-candidates with $s_0\le s-1$, which are constant on realization spaces by induction, and on pencils, where they are trivially true. If $[\mdr\le s-1]$ is false on $\cR(\Lambda)$, apply $(C_s)$.

For \GST: Theorem~\ref{thm:reduction} only involves $s_0\le s<K\le S+1$ when $d\le2S+4$. Choose $\cB$ by a rule depending only on the labelled lattice, e.g.\ the first $2s_0+3$ lines of $\cA_0$. Then for arrangements with isomorphic lattices the chosen sub-arrangements realize the same $s_0$-candidate.
\end{proof}

\begin{remark}\label{rem:galois}
If $\cR(\Lambda)$ is finite and consists of one orbit of $\operatorname{Gal}(\overline\QQ/\QQ)$, then $(C_s)$ holds for $\Lambda$: the dimensions $\dim_\CC D_0(\cA)_q$ are Galois invariant.
\end{remark}

%======================================================================
\section{The canonical foliation and the residual lemma}\label{sec:foliation}
%======================================================================

\subsection{Foliation and inflection polynomial}
A derivation $\theta\in D(\cB)_e\setminus SE$ defines a foliation of $\PP^2$. We call $\theta$ \emph{saturated} if the coefficients of $\theta\wedge E$ have no common factor, equivalently if its image in $T_{\PP^2}(e-1)$ has no divisorial zeros; in this case the foliation has degree $e$. A saturated $\theta$ also has isolated zeros as a section of $\mathcal T\langle\cB\rangle(e)\subseteq T_{\PP^2}(e-1)$, and we denote this zero scheme by $Z_\cB(\theta)$. The converse fails: for $\cB=\{xyz=0\}$, the derivation $\theta=-x\partial_x-y\partial_y+2z\partial_z$ is a nonzero element of $D_0(\cB)_1$ and has no zeros as a section of $\mathcal T\langle\cB\rangle(1)\cong\cO_{\PP^2}^2$, but $\theta\wedge E$ is divisible by $z$. The singular scheme of the foliation of a saturated $\theta$ is the zero scheme of its image in $T_{\PP^2}(e-1)$.

\begin{lemma}\label{lem:sat}
Let $\cA$ be critical with $s=\mdr(\cA)$ and $0\ne\theta\in D_0(\cA)_s$. Then $\theta$ is saturated and unique up to scalars and $SE$. No line of $\cA$ satisfies $\theta|_H\equiv0$, and
\[
\len Z_\cA(\theta)=c_2(\TA(s))=\nu(\cA).
\]
\end{lemma}
\begin{proof}
If $\theta|_H\equiv0$ for some $H$, Corollary~\ref{cor:gap} would give $r_H\ge d-s=s+3$, contradicting the candidate condition. If $\theta\wedge E$ had a common factor $h$, it would define a curve containing no line of $\cA$, and dividing by $h$ would produce a tangent derivation of smaller degree. Uniqueness: for another $\eta$ of degree $s$, $\det(E,\theta,\eta)$ has degree $2s+1<d$ and is divisible by $Q_\cA$, hence vanishes, so $\eta\wedge E=\lambda\,\theta\wedge E$ with $\lambda$ a constant by saturation. Finally $c_2(\TA(s))=c_2+sc_1+s^2=c_2-s(s+2)$.
\end{proof}

For $\theta=a\partial_x+b\partial_y+c\partial_z$ the \emph{inflection polynomial} is
\[
F_\theta=\det\begin{pmatrix}x&y&z\\ a&b&c\\ \theta(a)&\theta(b)&\theta(c)\end{pmatrix},\qquad\deg F_\theta=3e.
\]
It is unchanged by $\theta\mapsto\theta+hE$. Its zero locus is the inflection curve of the foliation \cite{Per01}. In the affine chart $z=1$, write the foliation as $X=A\partial_x+B\partial_y$ with $A=a-xc$, $B=b-yc$. Then
\begin{equation}\label{eq:affine}
F_\theta(x,y,1)=A\,X(B)-B\,X(A)=A^2B_x-B^2A_y+AB(B_y-A_x),
\end{equation}
as a direct computation shows.

\begin{proposition}\label{prop:residual}
If every line of $\cB$ is $\theta$-invariant, then $Q_\cB\mid F_\theta$. Moreover $F_\theta\ne0$ unless the foliation is a pencil of lines.
\end{proposition}
\begin{proof}
If $\theta(\alpha)=h\alpha$, then $\theta^2(\alpha)=\alpha(h^2+\theta(h))$, so along $\alpha=0$ the three rows of the determinant are linearly dependent. If $F_\theta\equiv0$, velocity and acceleration are everywhere collinear, so every leaf is contained in a line; two distinct leaf-lines meet only at singular points, which are finitely many, so the leaf-lines form a pencil.
\end{proof}

For a critical $\cA$ we write $F_\theta=Q_\cA R$ with $\deg R=3s-(2s+3)=s-3$; in particular $s\ge3$.

\subsection{Local frames}\label{ssec:frames}
Let $p$ be a point, $X$ the affine field of $\theta$ near $p$, and $Z=Z_\cB(\theta)$. The sheaf $T_{\PP^2}(-\log\cB)$ has the following local bases at $p$ (coordinates $x,y$ centred at $p$):
\begin{itemize}[leftmargin=*]
\item \emph{type O} ($p\notin\cB$): $\partial_x,\partial_y$; so $X=A\partial_x+B\partial_y$ and $I_{Z,p}=(A,B)$;
\item \emph{type S} ($p$ on exactly one line, $\{y=0\}$): $\partial_x,\ y\partial_y$; so $X=a\partial_x+b\,y\partial_y$ and $I_{Z,p}=(a,b)$;
\item \emph{type V} ($p$ a vertex of multiplicity $m\ge2$, local product $Q_p$ of degree $m$): $E_p=x\partial_x+y\partial_y$ and the Hamiltonian field $\eta_p=\partial_yQ_p\,\partial_x-\partial_xQ_p\,\partial_y$; so $X=aE_p+b\eta_p$ and $I_{Z,p}=(a,b)$.
\end{itemize}
In all cases $\det$ of the basis is a unit times the local equation of $\cB$, by Saito's criterion \cite{Sai80}. We write $D=F_\theta(x,y,1)$ and $n_p=\len Z_p$.

\begin{theorem}[residual lemma]\label{thm:residual}
Let $\theta\in D(\cB)_e$ be saturated, $Z=Z_\cB(\theta)$, and $F_\theta=Q_\cB R$. Then $R\in H^0(\cI_Z^2(3e-|\cB|))$. If $p\in Z$ is a vertex of $\cB$, then $\mult_pR\ge3$.
\end{theorem}
\begin{proof}
Let $\xi,\eta$ be a local basis as above, $X=f\xi+g\eta$, $I=(f,g)$. Let $\nabla$ be the flat connection of the affine chart, so that $X(X)=\nabla_XX$. For logarithmic $Y,W$ and an affine equation $\alpha$ of a line of $\cB$ we have $(\nabla_YW)(\alpha)=Y(W(\alpha))\in(\alpha)$, so $\nabla$ preserves logarithmic fields. Expanding,
\[
\nabla_XX=X(f)\xi+X(g)\eta+f(f\nabla_\xi\xi+g\nabla_\eta\xi)+g(f\nabla_\xi\eta+g\nabla_\eta\eta)=\varphi\xi+\psi\eta,
\]
and $\varphi,\psi\in I$ because $X(f)=f\xi(f)+g\eta(f)$, and similarly for $g$. Hence $D=\det(X,\nabla_XX)=(f\psi-g\varphi)\det(\xi,\eta)$ with $f\psi-g\varphi\in I^2$, and $R=D/Q_\cB$ is $(f\psi-g\varphi)$ times a unit.

At a vertex, $\xi$ and $\eta$ vanish at $p$. If $f,g\in\mm_p$, then $X(f),X(g)\in\mm_p^2$, and the remaining terms of $\varphi,\psi$ have coefficients in $\mm_p^2$. Hence $f\psi-g\varphi\in\mm_p^3$.
\end{proof}

\begin{corollary}\label{cor:first}
Let $\cA$ be critical and $Z=Z_\cA(\theta)$. Then $Z\subseteq\Sing R$. If $s\le4$, then $\cA$ is free. If $s=5$ and $Z\ne\varnothing$, then $R$ is a singular conic and $Z$ contains no vertex.
\end{corollary}
\begin{proof}
For $s=3,4$, $R$ has degree $0$ or $1$, so it does not lie in $\mm_p^2$; hence $Z=\varnothing$ and $\nu=0$. For $s=5$, a conic has no point of multiplicity $3$.
\end{proof}

\subsection{Initial forms of the residual}
In this subsection $\theta\in D(\cB)_e$ is saturated, $Z=Z_\cB(\theta)$ and $F_\theta=Q_\cB R$, as in Theorem~\ref{thm:residual}; a \emph{vertex} is a point of $\cB$ of multiplicity $\ge2$. The results are local, and they will be applied to $\cB=\cA$ and to $\cB=\cA\cup\{\ell\}$.

Let $p\in Z$. Since $R\in I_{Z,p}^2$, if $\mult_pR\le3$ then $I_{Z,p}\not\subseteq\mm_p^2$, so $Z_p$ is \emph{curvilinear}: $Z_p\cong\operatorname{Spec}\CC[t]/(t^{n_p})$, embedded along a smooth germ.

\begin{lemma}[vertex formula]\label{lem:vertexformula}
At a vertex $p$ of multiplicity $m$, with $X=aE_p+b\eta_p$, write $\operatorname{Hess}_{Q_p}(\eta_p,\eta_p)=Q_p\,k_p$. Then $D=Q_p\,R_p$ with
\[
R_p=m\bigl(b\,X(a)-a\,X(b)\bigr)-m(m-1)a^2b-k_pb^3 ,
\]
and $k_p$ is homogeneous of degree $2m-4$; for $m=2$ it is a nonzero constant. Near $p$ we have $Q_\cB=U\,Q_p$ with $U(p)\ne0$, so $R_p=U R$ is $R$ times a local unit.
\end{lemma}
\begin{proof}
Write $E=E_p$, $\eta=\eta_p$, $Q=Q_p$, $H=\operatorname{Hess}Q$. Then $\nabla_XX=X(a)E+X(b)\eta+aX+b\,(D\eta)X$, and $\det(E,\eta)=-mQ$ by Euler's formula. Since $D\eta$ is traceless, $\det(X,(D\eta)X)=-H(X,X)$. Moreover $H(E,E)=m(m-1)Q$, and $H(E,\eta)=(m-1)\nabla Q\cdot\eta=0$. Hence
\[
D=mQ\bigl(bX(a)-aX(b)\bigr)-m(m-1)Qa^2b-b^3H(\eta,\eta).
\]
Taking $a=0$, $b=1$ shows $Q\mid H(\eta,\eta)$. For two lines, $H(\eta,\eta)$ is a nonzero constant times $Q$.
\end{proof}

\begin{lemma}[initial forms]\label{lem:initial}
Let $p\in Z$ with $\mult_pR\le3$.
\begin{enumerate}[label=\textup{(L\arabic*)},leftmargin=3em]
\item If $\mult_pR=2$, then $n_p=1$ and $p$ is not a vertex.
\item If $p$ is of type O and $n_p\ge2$, then $\inn_3R$ is divisible by the tangent line of $Z_p$. If $n_p\ge3$, then $\inn_3R$ is the cube of that line.
\item If $p$ is of type S and $n_p\ge3$, then $\inn_3R$ is a cube.
\item If $p$ is a vertex of multiplicity $m\ge3$ and $n_p\ge2$, then the linear parts satisfy $b_1=\mu a_1$ with $\mu\ne0$, and $\inn_3R$ is a nonzero multiple of $\mu\,a_1^3$; precisely, $\inn_3R=-U(p)^{-1}m(m-1)\mu\,a_1^3$, where $Q_\cB=U\,Q_p$ near $p$. If $m=2$ and $n_p\ge2$, normalize $Q_p=l_1l_2$ so that the linear forms $l_1,l_2$ have determinant $1$; then either $b_1=\mu a_1$ with $\mu\notin\{0,\pm1\}$ and $\inn_3R=\varepsilon\,\mu(\mu^2-1)a_1^3$, or $a_1=0$ and $\inn_3R=\varepsilon'b_1^3$, where $\varepsilon,\varepsilon'$ are nonzero constants.
\end{enumerate}
\end{lemma}
\begin{proof}
(L1) At a vertex, Theorem~\ref{thm:residual} gives $\mult\ge3$. At a point of type O with linear part $J$ of $X$, $\inn_2D=\det(Jw,J^2w)=\det J\cdot\det(w,Jw)$. At a point of type S, with $X=uf\partial_u+g\partial_v$, $f=\alpha u+\beta v+\dots$ and $g=\gamma u+\delta v+\dots$, one computes
\[
\inn_2(D/u)=(\alpha\delta-\beta\gamma)\,u(\gamma u+\delta v).
\]
In both cases $n_p\ge2$ means that the determinant vanishes, and then $\mult_pR\ge3$.

(L2) Let $n_p\ge2$. Then $J$ has rank $1$, and we choose coordinates in which the tangent line of $Z_p$ is $\{x=0\}$; the linear parts of $A,B$ are then multiples of $x$.

If $A=x+A_2+\dots$ (after scaling) and, after a shear $y\mapsto y-cx$, $B=B_2+\dots$, then $\inn_3D=x(b_{20}x^2-b_{02}y^2)$. Here $n_p\ge3$ if and only if $b_{02}=0$.

If $A=A_2+\dots$ and $B=x+B_2+\dots$, then $\inn_3D=-x^2\partial_yA_2=-x^2(a_{11}x+2a_{02}y)$. Here $n_p\ge3$ if and only if $a_{02}=0$.

(L3) This is a direct computation of the same kind. In the frame $X=uf\partial_u+g\partial_v$, the condition $n_p\ge2$ means that the linear parts of $f,g$ are dependent. If $f_1=0$ and $g_1=\gamma u+\delta v$, then
\[
\inn_3(D/u)=-(\gamma u+\delta v)\bigl(\gamma f_{11}u^2-\delta f_{20}u^2+2\gamma f_{02}uv+\delta f_{02}v^2\bigr),
\]
and $n_p\ge3$ means $f_2(\delta,-\gamma)=0$, which makes the second factor proportional to $(\gamma u+\delta v)^2$. If $g_1=\mu f_1$, the analogous computation (Appendix~\ref{app:local}, item~(A0)) gives a cube as soon as $n_p\ge3$.

(L4) By Lemma~\ref{lem:vertexformula}, $R=U^{-1}R_p$ with $U(p)\ne0$, so $\inn_3R=U(p)^{-1}\inn_3R_p$; it suffices to compute $\inn_3R_p$. Apply Lemma~\ref{lem:vertexformula} with $a=a_1+\dots$, $b=b_1+\dots$, where $E_p(a_1)=a_1$ and $\eta_p$ has degree $m-1$.

For $m\ge3$, the terms of $m(bX(a)-aX(b))$ involving $\eta_p$ have order $\ge m+1\ge4$, so its cubic part is $m(a_1^2b_1-a_1^2b_1)=0$, the term $-m(m-1)a^2b$ contributes $-m(m-1)a_1^2b_1$, and $k_pb^3$ has order $\ge2m-1\ge5$. Moreover $n_p\ge2$ means that $a_1,b_1$ are dependent and, by curvilinearity, not both zero. If $a_1=0$, the cubic part vanishes, contradicting $\mult_pR\le3$.

For $m=2$, $\eta_p$ is linear and $k_p$ is constant, and
\[
\inn_3R_p=2\bigl(b_1^2\eta_p(a_1)-a_1b_1\eta_p(b_1)\bigr)-2a_1^2b_1-k_pb_1^3 .
\]
If $b_1=\mu a_1$, the first term vanishes and $\inn_3R_p=-\mu(2+k_p\mu^2)a_1^3$; if $a_1=0$, then $\inn_3R_p=-k_pb_1^3$. For $l_1,l_2$ of determinant $\Delta$ one has $k_p=-2\Delta^2$, so $\inn_3R_p=2\mu(\mu^2\Delta^2-1)a_1^3$, and the normalization $\Delta=1$ gives the two stated cases. In \S\ref{sec:fifteen} the two lines through a double point are normalized to $x,y$ or to $y,y-x$, both of determinant $\pm1$.
\end{proof}

\begin{lemma}\label{lem:twothirds}
Let $p\in Z$ and let $\ell$ be a line through $p$ with $\alpha_\ell^3\in I_{Z,p}^2$, $Z_p$ curvilinear. Then $\iota_p:=\len(Z_p\cap\ell)\ge\frac23n_p$, i.e.\ $n_p\le\lfloor\frac32\iota_p\rfloor$.
\end{lemma}
\begin{proof}
Write $I_{Z,p}=(y',t^{n})$ with $y'$ a local equation of the smooth germ carrying $Z_p$ and $t$ a parameter on it, $n=n_p$. Then $I^2=(y'^2,y't^n,t^{2n})$. Write $\alpha_\ell\equiv u_0(t)+y'u_1(t)\pmod{y'^2}$; then $\alpha_\ell^3\equiv u_0^3+3u_0^2u_1y'$. The condition $\alpha_\ell^3\in I^2$ gives $3\ord u_0\ge2n$, and $\iota_p=\min(n,\ord u_0)$.
\end{proof}

%======================================================================
\section{Defect at most one, and arrangements of at most fourteen lines}\label{sec:nu1}
%======================================================================

\begin{lemma}\label{lem:c1c2zero}
A rank-$2$ bundle $F$ on $\PP^2$ with $c_1=c_2=0$ and $H^0(F(-1))=0$ is trivial.
\end{lemma}
\begin{proof}
$F^*\cong F$, and Serre duality gives $H^2(F)=H^0(F(-3))^*=0$. Since $\chi(F)=2$, there is a section $\sigma\ne0$. It has no divisorial zeros, since $H^0(F(-1))=0$, and its zero scheme has length $c_2=0$. Hence $0\to\cO\to F\to\cO\to0$, which splits.
\end{proof}

\begin{proposition}\label{prop:nu1}
Let $\Lambda$ be an $s$-candidate all of whose realizations have $\mdr\ge s$, and suppose that $\Lambda$ has a critical realization $\cA$ with $\nu\le1$.
\begin{enumerate}[label=\textup{(\alph*)}]
\item If $\nu=0$, every realization of $\Lambda$ is critical.
\item If $\nu=1$, every realization of $\Lambda$ that is not critical is free with exponents $(s+1,s+1)$. Hence, if Terao's conjecture holds for the free arrangements of $2s+3$ lines with exponents $(s+1,s+1)$, every realization of $\Lambda$ is critical.
\end{enumerate}
\end{proposition}
\begin{proof}
Let $\cB\in\cR(\Lambda)$ with $\mdr(\cB)\ge s+1$, and put $F=\mathcal T\langle\cB\rangle(s+1)$. Then $c_1(F)=1-(2s+3)+2(s+1)=0$, $c_2(F)=c_2-(s+1)^2=s(s+2)+\nu-(s+1)^2=\nu-1$, and $H^0(F(-1))=D_0(\cB)_s=0$.

(a) If $\nu=0$, then $c_2(F)=-1$. Since $c_1(F)=0$ and $H^0(F(-1))=0$, $F$ is semistable, and Bogomolov's inequality gives $c_2(F)\ge0$: a contradiction. So no such $\cB$ exists.

(b) If $\nu=1$, then $c_2(F)=0$, and $\cB$ is free by Lemma~\ref{lem:c1c2zero}, with exponents $(s+1,s+1)$. If Terao's conjecture holds for $L(\cB)=\Lambda$, the critical realization $\cA$ is free too, contradicting $\nu=1$ (Lemma~\ref{lem:cand}(b)).
\end{proof}

\begin{theorem}[$13$ lines]\label{thm:nu13}
Every critical arrangement of $13$ lines has $\nu\le1$.
\end{theorem}
\begin{proof}
Here $s=5$, and we may assume $Z\ne\varnothing$. Then $R$ is a singular conic and $Z$ contains no vertex (Corollary~\ref{cor:first}).

\emph{$Z$ is reduced.} Let $p\in Z$. Since a conic has multiplicity $\le2$ at every point, (L1) gives $n_p=1$.

If $R=L_1L_2$ with $L_1\ne L_2$, then $\Sing R$ is a point, and $\nu\le1$. So let $R=L^2$ and $\ell=\{L=0\}$; then $Z\subset\ell$.

\emph{Case $\ell\in\cA$.} At $p\in Z$ take $\ell=\{u=0\}$, $X=uf\partial_u+g\partial_v$. Then $\inn_2(D/u)\propto u^2$ forces $\delta=0$ in the formula of (L1), so $g(0,v)=O(v^2)$ at $p$. Thus $t_\ell(\theta)\in H^0(\cO_\ell(6))$ vanishes at the $r_\ell$ vertices on $\ell$ and doubly at each point of $Z$. Hence $r_\ell+2\nu\le6$. Since $12=\sum_{p\in\ell}(m_p-1)\le5r_\ell$, we have $r_\ell\ge3$, and so $\nu\le1$.

\emph{Case $\ell\notin\cA$.} A point of $Z$ on a line $H\in\cA$ would be of type S, and the formula in (L1) would force $u\mid L^2$, i.e.\ $\ell=H$. So all points of $Z$ are of type O. At $p\in Z$, $\det(w,Jw)$ is a nonzero square, which forces $J=\left(\begin{smallmatrix}\lambda&\beta\\0&\lambda\end{smallmatrix}\right)$ with $\lambda\beta\ne0$ and eigenline $\ell=\{y=0\}$.

Write $X=A\partial_x+C\partial_y$. If $C(x,0)=cx^j+\dots$ with $c\ne0$ and $j\ge2$, then $D|_{y=0}=\lambda^2jc\,x^{j+1}+\dots$, contradicting $L^2\mid D$. Hence $\ell$ is invariant.

Choose an affine chart in which the point at infinity of $\ell$ is not in $Z$, and write $X=P\partial_x+yS\partial_y$, $p=P(x,0)$, $q=S(x,0)$. Here $p\not\equiv0$ by saturation, and $q(x_1)=\lambda\ne0$ at each point $(x_1,0)$ of $Z$. The coefficient of $y$ in $D$ is $p(pq'-qp'+q^2)$. It vanishes, so $(p/q)'=1$ and $p/q=x-x_0$. At each point of $Z$ we have $p=0$ and $q\ne0$, so $Z=\{(x_0,0)\}$.
\end{proof}

\begin{proof}[Proof of Theorem A]
By Theorem~\ref{thm:induction} with $S=5$ it suffices to prove $(C_s)$ for $s\le5$.
\begin{itemize}[leftmargin=*]
\item For $s\le2$ there are no critical arrangements. Indeed, by Lemma~\ref{lem:cand}(a) every line $H$ of an $s$-candidate satisfies $2s+2=\sum_{p\in H}(m_p-1)\le s\,r_H\le s(s+1)$. This is impossible for $s\le1$. For $s=2$ it forces every line to contain three triple points, so the only $2$-candidate is the Fano plane, which is not realizable over $\CC$.
\item For $s=3,4$, critical arrangements are free (Corollary~\ref{cor:first}), so Proposition~\ref{prop:nu1}(a) applies. No case of Terao's conjecture is needed.
\item For $s=5$, Theorem~\ref{thm:nu13} gives $\nu\le1$, and Proposition~\ref{prop:nu1} applies, using Terao's conjecture for free arrangements of $13$ lines with exponents $(6,6)$ \cite{DIM19}.
\end{itemize}
The statement about $14$ lines follows from the corollary to Proposition~\ref{prop:freecrit}.
\end{proof}

%======================================================================
\section{Fifteen lines: the defect is at most three}\label{sec:fifteen}
%======================================================================

Throughout this section $\cA$ is a critical $6$-candidate: $d=15$, $\mdr(\cA)=6$, and $r_H\le7$ for all $H$. By Lemma~\ref{lem:cand}(a), every multiplicity is $\le7$. Since $\sum_{p\in H}(m_p-1)=14$, this gives
\begin{equation}\label{eq:rH}
3\le r_H\le7\qquad\text{for every }H\in\cA .
\end{equation}
Let $\theta$, $\cF$ and $Z=Z_\cA(\theta)$ be as in \S\ref{sec:foliation}. Then $R$ is a cubic, $\nu=\len Z=43-X(\cA)$, and every point of $Z$ is a singular point of $R$ of multiplicity $\le3$, so $Z$ is curvilinear.

\begin{theorem}\label{thm:nu15}
Every critical $6$-candidate has $\nu\le3$, i.e.\ $X(\cA)\ge40$.
\end{theorem}

The proof occupies the rest of this section. The statements labelled \emph{computer-assisted} are verified by finitely many exact computations, listed with their outcomes in Appendix~\ref{app:local}.

\subsection{Invariant lines: a Riccati equation}\label{ssec:riccati}
Let $\ell$ be a $\theta$-invariant line, not necessarily in $\cA$. In an affine chart with $\ell=\{u=0\}$ write $X=uf\partial_u+g\partial_v$, $f=\sum u^if_i(v)$, $g=\sum u^ig_i(v)$. Here $g_0\ne0$ by saturation, and $t_\ell(\theta)=g_0\partial_v$. By \eqref{eq:affine}, $D=u\widetilde D$ with
\begin{equation}\label{eq:Dtilde}
\widetilde D=uf^2g_u+fgg_v-f^2g-ufgf_u-g^2f_v .
\end{equation}

\begin{lemma}[Riccati]\label{lem:riccati}
Suppose $u\mid\widetilde D$.
\begin{enumerate}[label=\textup{(\alph*)}]
\item There is a unique point $p_0\in\ell$ such that the defining $1$-forms of $\cF$ and of the pencil of lines through $p_0$ agree to first order along $\ell$, i.e.\ their wedge product vanishes modulo $u^2$. The point $p_0$ is a singular point of $\cF$. In an affine chart whose line at infinity passes through $p_0$, one has $f_0\equiv0$.
\item In such a chart, if moreover $u^2\mid\widetilde D$, then $f_1=\kappa g_0$ for a constant $\kappa$. If $u^3\mid\widetilde D$, then also $f_2=\kappa g_1-(2\kappa^2v+\kappa')g_0$ for a constant $\kappa'$.
\item In such a chart, let $j\in\{1,2,3\}$ and suppose $u^j\mid\widetilde D$. Then the coefficient of $u^j$ in $\widetilde D$ is divisible by $g_0$ in $\CC[v]$.
\end{enumerate}
\end{lemma}
\begin{proof}
The coefficient of $u^0$ in $\widetilde D$ is $g_0(f_0g_0'-g_0f_0'-f_0^2)$. So either $f_0\equiv0$, or $(g_0/f_0)'=1$ and $g_0=(v-v_0)f_0$.

The foliation is defined by $\omega=g\,du-uf\,dv$. The pencil through $(0,v_0)$ is defined by $(v-v_0)du-u\,dv$, and
\[
\omega\wedge\bigl((v-v_0)du-u\,dv\bigr)=u\bigl((v-v_0)f-g\bigr)\,du\wedge dv,
\]
which vanishes modulo $u^2$ if and only if $g_0=(v-v_0)f_0$. The pencil through the point at infinity of $\ell$ is defined by $du$, and $\omega\wedge du=uf\,du\wedge dv$, which vanishes modulo $u^2$ if and only if $f_0\equiv0$. Since $g_0\ne0$, at most one point satisfies the condition, and the Riccati equation shows that one does.

The condition is intrinsic. In homogeneous terms it says that $\omega\wedge\omega_{p}$ vanishes to order $\ge2$ along $\ell$, where $\omega$ and $\omega_p$ are the twisted $1$-forms defining $\cF$ and the pencil through $p$. Both forms change by units under a change of affine chart, and the vanishing order along $\ell$ can be read off on any dense open subset of $\ell$. So the condition is preserved by projective transformations preserving $\ell$, and the two computations above cover every point of $\ell$. In a chart where $p_0$ is at infinity we therefore get $f_0\equiv0$.

That $p_0$ is singular: if $g_0=(v-v_0)f_0$, then $g_0(v_0)=0$. If $f_0\equiv0$, recall that $s=6$ in this section and write the degree-$7$ part of the affine field as $-c(u,v,0)\cdot(u,v)$, with $c$ the $\partial_z$-coefficient of $\theta$ (in the coordinates of the chart); then the $v^6$-coefficient of $f_0$ and the $v^7$-coefficient of $g_0$ both equal $-c(0,1,0)$. Hence $\deg g_0\le6$, i.e.\ $t_\ell$ vanishes at infinity.

(b) and (c) follow from the expansion of \eqref{eq:Dtilde}, whose coefficients are, with $f_0=0$:
\begin{align*}
[u^1]&=g_0\,(f_1g_0'-g_0f_1'),\\
[u^2]&=g_0\bigl(-2\kappa^2g_0^2+\kappa(g_0g_1'-g_1g_0')+(f_2g_0'-g_0f_2')\bigr)\quad\text{if }f_1=\kappa g_0,
\end{align*}
together with the analogous divisibility of $[u^3]$ by $g_0$ when $f_1,f_2$ are as stated. These identities are verified symbolically, with $f_i,g_i$ arbitrary functions (Appendix~\ref{app:local}, item~(A1)). The vanishing of $[u^1]$ gives $(f_1/g_0)'=0$. The vanishing of $[u^2]$ gives $(\kappa g_1-f_2)/g_0=2\kappa^2v+\kappa'$.
\end{proof}

\begin{corollary}\label{cor:riccati}
In the situation of Lemma~\ref{lem:riccati}(c), let $\widetilde D=u^j\Phi$ with $\Phi(0,v)\not\equiv0$, $j\in\{1,2,3\}$. Let $q\in\ell$ be a zero of $t_\ell(\theta)$ with $q\ne p_0$. Then
\[
1\le\ord_qt_\ell(\theta)\le\ord_q\Phi|_\ell .
\]
\end{corollary}
\begin{proof}
Use the chart with $p_0$ at infinity. Then $q$ is a root of $g_0$, and $g_0$ divides $[u^j]\widetilde D=\Phi(0,v)$.
\end{proof}

\begin{lemma}\label{lem:tangentorder}
Let $\ell$ be invariant and $q\in\ell\cap Z$, with $\iota_q=\len(Z_q\cap\ell)$.
\begin{enumerate}[label=\textup{(\alph*)}]
\item If $\ell\in\cA$ and $q$ is a vertex of multiplicity $m$ with $X=aE_q+b\eta_q$, then, in coordinates centred at $q$ with $\ell=\{x=0\}$, the other lines through $q$ given by $y=t_ix$, $a_0=a|_\ell$ and $b_0=b|_\ell$,
\[
t_\ell=\bigl(a_0y-b_0y^{m-1}\bigr)\partial_y,\qquad f_0=a_0+(m-1)y^{m-2}b_0,
\]
where $X_x=xf$. In particular $\ord_qt_\ell\ge1+\iota_q$.
\item If $\ell\notin\cA$ and $q$ is a point of $\cA$, then $\ord_qt_\ell=1+\iota_q$.
\end{enumerate}
\end{lemma}
\begin{proof}
(a) Write $Q_q=x\,\Pi$ with $\ell=\{x=0\}$ and $\Pi=\prod(y-t_ix)$. Then $X_x=x(a+b\Pi_y)$ and $X_y=ay-b(\Pi+x\Pi_x)$; restrict to $x=0$. Moreover $Z_q\cap\ell=V(a_0,b_0)$.

(b) If $q$ is a vertex of $\cA$, invariance of $\ell$ forces $b_0=0$. Then $t_\ell=a_0y$ and $Z_q\cap\ell=V(a_0)$. If $q$ lies on one line $\{y=0\}$ of $\cA$, then $X=a\partial_x+b\,y\partial_y$ with $a=x\tilde a$, $t_\ell=yb_0$, and $Z_q\cap\ell=V(b_0)$.
\end{proof}

\subsection{Residual not a cube}\label{ssec:notcube}

\begin{proposition}\label{prop:notcube}
If $R$ is not the cube of a linear form, then $\nu\le3$.
\end{proposition}
\begin{proof}
\emph{$R$ reduced.} At a point of multiplicity $2$, (L1) gives $n_p=1$. If $R$ has no triple point, it has at most three singular points, so $\nu\le3$. If $R$ is three concurrent lines, then $\operatorname{Supp}Z\subseteq\{p\}$, and $\inn_3R$ is a product of three distinct lines. By (L2) and (L3), $n_p\le2$ at points of types O and S. By (L4), $n_p\le1$ at a vertex.

\emph{$R=L_1L_2^2$, $L_1\ne L_2$.} Then $\operatorname{Supp}Z\subseteq\ell_2=\{L_2=0\}$. At $p_0'=\ell_1\cap\ell_2$, the cubic $\inn_3R$ is not a cube, so $n\le2$ there, as above. At the other points of $\ell_2$, $R$ has multiplicity $2$, so $Z$ is reduced there and avoids vertices. We show that $Z\cap(\ell_2\setminus p_0')$ has at most one point, which gives $\nu\le3$.
\begin{itemize}[leftmargin=*]
\item \emph{$\ell_2\in\cA$.} Here $\widetilde D=u^2\Phi$ with $\Phi=(Q/u)L_1$ up to a constant. By Corollary~\ref{cor:riccati}, the zeros of $t_{\ell_2}$ other than $p_0$ are vertices or $p_0'$. A point $q\in\operatorname{Supp}Z\setminus\{p_0'\}$ of type S on $\ell_2$ is a zero of $t_{\ell_2}$, so $q=p_0$.
\item \emph{$\ell_2\notin\cA$ invariant.} Here $\widetilde D=u\Phi$ with $\Phi=QL_1$. Points of $Z$ of type O on $\ell_2$ are zeros of $t_{\ell_2}$ not in $\cA$, so they equal $p_0$ or $p_0'$. A point $q\in Z\setminus\{p_0,p_0'\}$ of type S on $\ell_2$ would satisfy $\ord_qt_{\ell_2}=1+\iota_q\ge2$ by Lemma~\ref{lem:tangentorder}(b), whereas $\ord_q\Phi|_{\ell_2}=1$ since $q\ne p_0'$, contradicting Corollary~\ref{cor:riccati}, which applies since $q\ne p_0$. So $Z\cap(\ell_2\setminus p_0')\subseteq\{p_0\}$.
\item \emph{$\ell_2$ not invariant.} We show $Z\cap(\ell_2\setminus p_0')=\varnothing$.

Let $q$ be such a point, $\ell_2=\{x=0\}$. If $q$ is of type O, then $\inn_2D=\det J\det(w,Jw)\propto x^2$ forces $J=\left(\begin{smallmatrix}\lambda&0\\ \beta&\lambda\end{smallmatrix}\right)$ with $\lambda\beta\ne0$. Then $A(0,y)=\alpha y^j+\dots$ with $j\ge2$ and $B(0,y)=\lambda y+\dots$. The coefficient of $x^0$ in $D$ starts with $-j\alpha\lambda^2y^{j+1}$. Since $x^2\mid D$, we get $A(0,y)\equiv0$, i.e.\ $\ell_2$ is invariant: a contradiction.

If $q$ is of type S, the coefficient of $x^2y$ in $D$ vanishes identically when $q\in Z$ (Appendix~\ref{app:local}, item~(C4)). But $D=y\,x^2\cdot(\text{unit})$ near $q$.
\end{itemize}
\end{proof}

\subsection{\texorpdfstring{$R=L^3$}{R = L\^3}: generalities}
From now on $R=L^3$, $\ell=\{L=0\}$. Then $\operatorname{Supp}Z\subseteq\ell$, every point of $Z$ is curvilinear, and Lemma~\ref{lem:twothirds} gives
\begin{equation}\label{eq:twothirds}
n_p\le\bigl\lfloor\tfrac32\iota_p\bigr\rfloor,\qquad \iota_p=\len(Z_p\cap\ell)\ge1\quad(p\in Z).
\end{equation}
At $p\in\ell$, $D=Q_p\,x^3\,U$ with $U(p)\ne0$, where $x=0$ is the local equation of $\ell$. In particular $\inn_3R=\varepsilon\,x^3$ with $\varepsilon\ne0$. At a vertex $p\in Z$ with $n_p\ge2$, (L4) therefore gives $a_1=\alpha x$, and after scaling $X$ we may write
\begin{equation}\label{eq:branches}
\begin{aligned}
&\text{($\mu$-branch)} && a=x+O(2),\quad b=\mu a+c,\quad c\in\mm_p^2,\quad \mu\ne0;\\
&\text{(B-branch, only for $m=2$)} && a\in\mm_p^2,\quad b=x+c,\quad c\in\mm_p^2 .
\end{aligned}
\end{equation}
We use coordinates centred at $p$ with $\ell=\{x=0\}$. The other lines through $p$ are $y=t_ix$, normalized by $y\mapsto\beta y+\gamma x$ so that two of the slopes are $0$ and $1$; the remaining slopes are kept as parameters. Pure coefficients are those of $a_0=a(0,y)$, $c_0=c(0,y)$, i.e.\ of $x^0y^j$.

\subsection{Case (a): \texorpdfstring{$\ell=H\in\cA$}{l = H in A}}\label{ssec:casea}
Here $D=Q\,\lambda u^3$ in a chart with $H=\{u=0\}$, so $\widetilde D=u^3\Phi$ with $\Phi=\lambda Q/u$. Lemma~\ref{lem:riccati} applies with $j=3$. Fix a chart with $p_0$ at infinity. By Corollary~\ref{cor:riccati}, the zeros of $t_H\in H^0(\cO_H(7))$ other than $p_0$ are vertices $q$, with
\begin{equation}\label{eq:aorders}
1+\iota_q\le\ord_qt_H\le m_q-1\qquad(q\neq p_0),
\end{equation}
where the left inequality is Lemma~\ref{lem:tangentorder}(a) and $\iota_q:=0$ if $q\notin Z$. In particular no point of type S other than $p_0$ is singular, and double points other than $p_0$ are not in $Z$.

\begin{lemma}[computer-assisted]\label{lem:casea}
\begin{enumerate}[label=\textup{(\alph*)}]
\item $p_0\notin Z$ if $p_0$ is a point of type S.
\item A vertex $q\ne p_0$ in $Z$ has $n_q=1$.
\item If $p_0$ is a vertex in $Z$ with $n_{p_0}\ge2$: if $m_{p_0}\ge3$, then $\iota_{p_0}\ge5$; if $m_{p_0}=2$, then $p_0$ is in the $\mu$-branch and $n_{p_0}=\iota_{p_0}=\ord_{p_0}t_H-1$.
\end{enumerate}
\end{lemma}
\begin{proof}
(a) Let $p_0$ be of type S, in a chart where $p_0=0$ is finite, and write $X=ax\,\partial_x+b\,\partial_y$, $I_{Z,p_0}=(a,b)$. Put $\Delta=a_{01}b_{10}-a_{10}b_{01}$. The coefficients of $x^2y$ and $x^3$ in $D$ are $-b_{01}\Delta$ and $-b_{10}\Delta$, and the coefficient of $x^4$ vanishes when $b_{10}=b_{01}=0$ (item~(A2)). Since $D=x^4\cdot(\text{unit})$, we get $\Delta=0$ and $(b_{10},b_{01})\ne(0,0)$. Hence $n_{p_0}\ge2$, and by Lemma~\ref{lem:twothirds} $Z_{p_0}$ is tangent to $H$. So $b_{01}=a_{01}=0$, and after scaling $b=x+c$, $a=\alpha x+O(2)$ with $c\in\mm^2$. The Riccati relation $g_0=yf_0$ reads $b_0=ya_0$, i.e.\ $c_0=ya_0$. Induction on $j$: if $a_{0,i}=0$ for $i<j$, then $c_{0,j}=a_{0,j-1}=0$, and the coefficient of $x^3y^{j-1}$ in $D$ equals $-j\,a_{0,j}$ (item~(A2)). Hence $a_0\equiv b_0\equiv0$ up to order $7$, so $\ord_{p_0}t_H\ge8$. This is impossible.

(b) By \eqref{eq:aorders}, $\iota_q\le m_q-2$, so double points are excluded. Let $m=m_q\ge3$ and $n_q\ge2$; we are in the $\mu$-branch. In the chart with $p_0$ at infinity $f_0\equiv0$, i.e.\ $a_0=-(m-1)y^{m-2}b_0$ by Lemma~\ref{lem:tangentorder}(a). Since $b_0=\mu a_0+c_0$, this expresses $a_0$ through $c_0$, and $\ord a_0>\ord c_0=\ord b_0$, so $\iota_q=\ord c_0\ge2$. For $m=3$ this contradicts $\iota_q\le1$. For $4\le m\le7$, induction: if $c_{0,i}=0$ for $i<j\le m-2$, the coefficient of $x^3y^{j+m-1}$ in $D$ is $-m(m+j-2)\,c_{0,j}$ (item~(A3)). Hence $\iota_q\ge m-1$, contradicting \eqref{eq:aorders}.

(c) In a chart where $p_0$ is finite, $g_0=yf_0$; by Lemma~\ref{lem:tangentorder}(a) this reads $m\,y^{m-1}b_0=0$, so $b_0\equiv0$. Then $t_H=a_0y$ and $\iota_{p_0}=\ord a_0$.

For $m\ge3$ we are in the $\mu$-branch, and by induction the coefficient of $x^3y^{j+m-1}$ is $-m(2m-1-j)\mu\,a_{0,j}$ for $2\le j\le4$ (item~(A4)). Hence $\ord a_0\ge5$.

For $m=2$ in the B-branch, the coefficient of $x^3y^{j+1}$ is $-2(j+1)a_{0,j}$ for $2\le j\le7$ (item~(A5)), so $\ord_{p_0}t_H\ge9>7$: impossible. In the $\mu$-branch, $c_0=-\mu a_0$. On the smooth curve $a=0$ we have $x=O(y^{\iota})$ with $\iota=\ord a_0=\iota_{p_0}$, and $c|_{a=0}=-\mu a_0+O(y^{\iota+1})$, since $\partial_xc(0,y)\in(y)$. Hence $n_{p_0}=\iota_{p_0}$.
\end{proof}

\begin{proposition}\label{prop:casea}
In case (a), $\nu\le3$.
\end{proposition}
\begin{proof}
Let $V$ be the set of vertices on $H$, $r_H=|V|\ge3$, and $Z'=Z\setminus\{p_0\}$, so that $Z'\subseteq V$. We have $\sum_{q\in H}\ord_qt_H=7$, $\ord_qt_H\ge1$ for $q\in V\cup\{p_0\}$, and $\ord_qt_H\ge2$ for $q\in Z'$, where $n_q=1$.
\begin{itemize}[leftmargin=*]
\item If $p_0\notin Z$, put $\rho=|V\setminus\{p_0\}|$. Counting orders gives $7\ge1+\rho+|Z'|$, and $|Z'|\le\rho$. Hence $\nu=|Z'|\le\min(\rho,6-\rho)\le3$.
\item If $p_0\in Z$ with $n_{p_0}=1$, then $p_0$ is a vertex by Lemma~\ref{lem:casea}(a), $\ord_{p_0}t_H\ge2$, and $|Z'|\le\min(r_H-1,\,7-2-(r_H-1))\le2$. Hence $\nu=1+|Z'|\le3$.
\item If $n_{p_0}\ge2$ and $m_{p_0}\ge3$, then $\ord_{p_0}t_H\ge6$, leaving at most $1$ for the $r_H-1\ge2$ other vertices: impossible.
\item If $n_{p_0}\ge2$ and $m_{p_0}=2$, the other vertices carry $\sum(m_q-1)=13$ with $m_q-1\le6$, so $r_H-1\ge3$. Then $(n_{p_0}+1)+(r_H-1)+|Z'|\le7$, and $\nu=n_{p_0}+|Z'|\le7-r_H\le3$.
\end{itemize}
\end{proof}

\subsection{Case (b): \texorpdfstring{$\ell\notin\cA$}{l not in A}, invariant}\label{ssec:caseb}
Let $r=r_\ell$ be the number of points of $\ell$ on lines of $\cA$, and $m_p$ the number of lines of $\cA$ through $p\in\ell$. Then $\sum_pm_p=15$, $m_p\le7$, hence $r\ge3$. Put $\cA'=\cA\cup\{\ell\}$. Then $\theta\in D(\cA')_6$; it is saturated, since saturation only involves $\theta\wedge E$, and $F_\theta=Q_{\cA'}L^2$.

\begin{lemma}\label{lem:Zprime}
Let $Z'=Z_{\cA'}(\theta)$. Then $\len Z'=\nu+r-7$ and $\len Z'\le1$.
\end{lemma}
\begin{proof}
For $\cA'$ we have $\tau(\cA')=\tau(\cA)+\sum_{p\in\ell}(2m_p-1)=\tau(\cA)+30-r$. Hence $c_2(\mathcal T\langle\cA'\rangle(6))=225-\tau(\cA')-90+36=\nu+r-7$, using $\tau(\cA)=148-\nu$.

By Theorem~\ref{thm:residual} applied to $\cA'$, we have $L^2\in I_{Z'}^2$, so $Z'\subset\ell$. Since $\mult L^2=2$, (L1) shows that $Z'$ is reduced and avoids the vertices of $\cA'$, i.e.\ all of $\ell\cap\cA$. Hence $Z'$ consists of singular points of $\cF$ on $\ell\setminus\cA$.

On the other hand, $D=Q\lambda u^3$ gives $\widetilde D=u^2\Phi$ with $\Phi=\lambda Q$. By Corollary~\ref{cor:riccati}, every zero of $t_\ell$ other than $p_0$ lies in $\cA$. So $Z'\subseteq\{p_0\}$.
\end{proof}

\begin{lemma}[Ziegler restriction]\label{lem:ziegler}
Consider the multiarrangement on $\ell\cong\PP^1$ given by the points $p\in\ell\cap\cA$ with multiplicities $m_p$. It admits a nonzero derivation of degree $6$.
\end{lemma}
\begin{proof}
Write $\theta(\alpha_\ell)=\alpha_\ell\varphi$ and put $\tilde\theta=\theta-\varphi E\in D(\cA')$, so that $\tilde\theta(\alpha_\ell)=0$. By Ziegler \cite{Zie89}, the restriction of $\tilde\theta$ to $\ell$ lies in the module of the multirestriction $(\cA'^{\ell},m^\ell)$, whose multiplicities are $|\cA'_p|-1=m_p$. If this restriction were $0$, then $\tilde\theta=\alpha_\ell\theta''$ with $\theta''\in D(\cA)_5=S_4E$, since $\mdr(\cA)=6$. Then $\theta\in SE$, a contradiction.
\end{proof}

\begin{lemma}\label{lem:multiarr}
If $r\in\{3,4\}$, the multiarrangement of Lemma~\ref{lem:ziegler} has no nonzero derivation of degree $6$.
\end{lemma}
\begin{proof}
It is the Ziegler restriction of $\cA'$ onto $\ell$, with $|m|=15\ge2m_p$ for every $p$, since $m_p\le7$. By a theorem of Abe \cite[Thm.~3.5]{Abe13} (see also \cite[Thm.~2.8]{Abe14}), its exponents $d_1\le d_2$ satisfy $d_2-d_1\le r-2\le2$. Since $d_1+d_2=15$ is odd, $(d_1,d_2)=(7,8)$. For $r=3$ this is also Wakamiko's formula \cite{Wak07}.

As an independent check for $r=4$, put the points at $0,\infty,1,\lambda$, with the two largest multiplicities at $0,\infty$, and write $p=y^{m_1}P_1$, $q=z^{m_2}Q_1$. The conditions at $1$ and $\lambda$ form a system of $m_3+m_4$ linear equations in $14-m_1-m_2$ unknowns, one more equation than unknowns. For each of the $16$ multiplicity types, the greatest common divisor of the maximal minors, a polynomial in $\lambda$, has no roots other than $0$ and $1$ (item~(B1)).
\end{proof}

\begin{proposition}\label{prop:caseb}
In case (b), $\nu\le3$.
\end{proposition}
\begin{proof}
By Lemmas~\ref{lem:ziegler} and~\ref{lem:multiarr}, $r\ge5$. By Lemma~\ref{lem:Zprime}, $\nu=7-r+\len Z'\le8-r\le3$.
\end{proof}

\subsection{Case (c): \texorpdfstring{$\ell$}{l} not invariant}\label{ssec:casec}
Let $r$ be the number of points of $\ell$ on lines of $\cA$, and $\sigma_\ell=\sum_{p\in\ell}(m_p-1)=15-r$, the sum over $p\in\ell\cap\cA$. The normal restriction $n_\ell(\theta)\in H^0(\cO_\ell(6))$ is nonzero, and by Lemma~\ref{lem:restr}(b) it vanishes to order $\ge m_p-1$ at each $p\in\ell\cap\cA$. Define
\[
\beta_p=\ord_pn_\ell-(m_p-1)\ \ (p\in\cA),\qquad\beta_p=\ord_pn_\ell\ \ (p\notin\cA).
\]
Then
\begin{equation}\label{eq:budget}
\sum_{p\in\ell}\beta_p=6-\sigma_\ell=r-9\ \ge0 .
\end{equation}
For $p\in Z$ we have $\beta_p\ge\iota_p\ge1$. Indeed, at a vertex $n_\ell=b_0\cdot\eta_p(\alpha_\ell)|_\ell$ and $\eta_p(\alpha_\ell)|_\ell$ has order exactly $m_p-1$, since $\ell\notin\cA$; at points of types S and O, $n_\ell=a_0$, resp.\ $A(0,y)$.

\begin{lemma}\label{lem:casecO}
$Z$ contains no point of type O.
\end{lemma}
\begin{proof}
Let $p\in Z$ be of type O, $\ell=\{x=0\}$, and $X=A\partial_x+B\partial_y$. The conditions $x^3\mid D$ say that the coefficients $D_{(0)},D_{(1)},D_{(2)}$ of $x^0,x^1,x^2$ in $D$ vanish identically in $y$. We write $A_k(y),B_k(y)$ for the coefficients of $x^k$; in particular $A_0=n_\ell$ locally, and $A_0\not\equiv0$ as $\ell$ is not invariant. We use
\begin{align*}
D_{(0)}&=A_0^2B_1-A_0A_1B_0+A_0B_0B_0'-B_0^2A_0',\\
D_{(1)}&=2A_0^2B_2+A_0A_1B_1-2A_0A_2B_0+A_0B_0B_1'+A_0B_1B_0'\\
&\quad-A_1^2B_0+A_1B_0B_0'-B_0^2A_1'-2B_0B_1A_0',\\
D_{(2)}&=3A_0^2B_3+3A_0A_1B_2-3A_0A_3B_0+A_0B_0B_2'+A_0B_1B_1'+A_0B_2B_0'\\
&\quad-3A_1A_2B_0+A_1B_0B_1'+A_1B_1B_0'+A_2B_0B_0'-B_0^2A_2'-2B_0B_1A_1'\\
&\quad-2B_0B_2A_0'-B_1^2A_0'.
\end{align*}
Let $\varpi=\ord A_0<\infty$ and $\varsigma=\ord B_0$, with $\varsigma=\infty$ if $B_0\equiv0$.

If $n_p=1$, the linear part $J$ of $X$ is invertible, and $\inn_2D=\det J\cdot\det(w,Jw)=0$ (see the proof of (L1)) forces $J=\lambda I$ with $\lambda\ne0$. So $\varsigma=1$, $A_1(0)=\lambda$ and $\varpi\ge2$. The coefficient of $y^{\varpi+1}$ in $D_{(0)}$ is $-\varpi\,a\lambda^2$, where $A_0=ay^\varpi+\dots$; hence $a=0$, a contradiction.

If $n_p\ge2$, then (L2) and $\inn_3D=\varepsilon x^3$ show that the tangent line of $Z_p$ is $\ell$, and the linear parts of $A,B$ are multiples of $x$.

\emph{Saddle-node} ($A=x+\dots$, $B=O(2)$, after a shear). Here $\inn_3D=x(b_{20}x^2-b_{02}y^2)$, so $b_{02}=0$ and $b_{20}\ne0$; hence $\varsigma\ge3$. In $D_{(0)}$, the term $A_0A_1B_0$ has exact order $\varpi+\varsigma$, while $A_0^2B_1$ has order $\ge2\varpi+1$ and $B_0(A_0B_0'-B_0A_0')$ has order $\ge\varpi+2\varsigma-1$. Hence $\varsigma\ge\varpi+1$. Since moreover $B_1\in(y)$, $A_1(0)=1$ and $B_2(0)=b_{20}\ne0$, every term of $D_{(2)}$ has order $\ge\varpi+1$ except $3A_0A_1B_2$, which has exact order $\varpi$. This is a contradiction.

\emph{Nilpotent} ($A\in\mm^2$, $B=x+\dots$). Here $\inn_3D=-x^2\partial_yA_2$ forces $A_2=a_{20}x^2-cxy$ with $c\ne0$. So $A_1=-cy+\dots$, $B_1=1+\dots$, $\varpi\ge3$ and $\varsigma\ge2$. Moreover $\varsigma<\infty$: if $B_0\equiv0$, then $D_{(0)}=A_0^2B_1\not\equiv0$.
\begin{itemize}[leftmargin=*]
\item The $y^2$-coefficient of $D_{(2)}$ is $-3a_{03}$, so $\varpi\ge4$.
\item If $\varsigma=2$, write $B_0=by^2+\dots$. The $y^4$-coefficient of $D_{(1)}$ is $-cb(c+b)$, so $b=-c$. The $y^{\varpi+3}$-coefficient of $D_{(0)}$ is $ab\bigl(c+(2-\varpi)b\bigr)=abc(\varpi-1)\ne0$: a contradiction.
\item If $\varsigma\ge3$, comparing orders in $D_{(0)}$ gives $\varpi=\varsigma+1$ and $a=-cb$. The $y^\varsigma$-coefficient of $D_{(2)}$ is then $-(\varsigma+1)a+(2-\varsigma)cb=3cb\ne0$: a contradiction.
\end{itemize}
The expansions of the saddle-node case are checked in item~(C1) of Appendix~\ref{app:local}; the nilpotent case is the explicit computation above.
\end{proof}

\begin{lemma}[computer-assisted]\label{lem:casec}
\begin{enumerate}[label=\textup{(\alph*)}]
\item $Z$ contains no point of type S.
\item A vertex $p\in Z$ with $m_p=2$ has $n_p=1$, unless $\beta_p\ge8$.
\item A vertex $p\in Z$ with $m_p=3$ has $n_p=1$, unless $\beta_p\ge6$.
\item A vertex $p\in Z$ with $m_p=4$ and $n_p\ge2$ has $\beta_p\ge3$, and $n_p=3$ if $\beta_p=3$.
\end{enumerate}
\end{lemma}
\begin{proof}
(a) Let $X=a\partial_x+b\,y\partial_y$ with $p\in Z$. The coefficients of $xy^2$ and $x^2y^2$ in $D$ and the nondegeneracy (coefficient of $x^3y\ne0$) force, after scaling, $a=x+O(2)$ and $b=\mu a+c$ with $c\in\mm^2$ and $c_{20}+\mu^2\ne0$ (item~(C2)). By induction, the coefficients of $xy^{j+1}$ and $x^2y^{j+1}$ are $-c_{0,j}$ and $3(c_{20}+\mu^2)a_{0,j}$ once the lower pure coefficients vanish, for $2\le j\le7$. Hence $n_\ell=a_0$ vanishes to order $\ge8>6$, so $n_\ell=0$: a contradiction.

(b) Both branches of \eqref{eq:branches}: by induction all pure coefficients up to order $7$ vanish (item~(C3)). Hence $b_0\in(y^8)$.

(c) $\mu$-branch, slopes $0,1,s_3$. The coefficients of $x^2y^5$, $xy^7$ and $y^9$ give $c_{0,2}=-2\mu a_{0,2}$, $a_{0,2}(2a_{0,2}-3\mu)=0$ and $a_{0,2}^2(a_{0,2}+3\mu)=0$. Hence $a_{0,2}=c_{0,2}=0$. For $3\le j\le5$, the coefficient of $x^2y^{j+3}$ expresses $c_{0,j}$ as a multiple of $\mu a_{0,j}$, and then the coefficient of $xy^{j+5}$ is a nonzero multiple of $\mu^2a_{0,j}$ (item~(C3)). Hence $b_0=\mu a_0+c_0\in(y^6)$.

(d) $\mu$-branch, slopes $0,1,s_3,s_4$. Here $c_{0,2}=-\frac94\mu a_{0,2}$ and $27\mu a_{0,2}^2=0$, so $a_{0,2}=c_{0,2}=0$. Then $c_{0,3}=-\frac95\mu a_{0,3}$, so $b_0=-\frac45\mu a_{0,3}y^3+O(y^4)$. If $\beta_p=3$, then $a_{0,3}\ne0$. On $a=0$ we have $x=-a_{0,3}y^3+\dots$, so $c|_{a=0}=c_{0,3}y^3+\dots$, and hence $n_p=3$.
\end{proof}

\begin{proposition}\label{prop:casec}
In case (c), $\nu\le3$.
\end{proposition}
\begin{proof}
By Lemmas~\ref{lem:casecO} and~\ref{lem:casec}(a), $Z$ consists of vertices on $\ell$. Their $\beta_p\ge1$ satisfy $\sum\beta_p\le r-9=6-\sigma_\ell$ by \eqref{eq:budget}.
\begin{itemize}[leftmargin=*]
\item If no vertex of multiplicity $\ge4$ lies in $Z$, Lemma~\ref{lem:casec}(b),(c) gives $n_p=1$; note that $\beta_p\le6-\sigma_\ell\le5$, resp.\ $\le4$. Hence $\nu\le\min(\sigma_\ell,6-\sigma_\ell)\le3$.
\item If $p\in Z$ has $m_p=4$, then $\sigma_\ell\ge3$ and $\sum\beta\le3$. Either $n_p=1$ and $\nu\le1+2$, or $\beta_p=3$, which exhausts the budget, and then $\nu=n_p=3$.
\item If $p\in Z$ has $m_p\ge5$, then $\sum\beta\le2$. Either a single point of $Z$ has $\beta=2$ and $n\le3$ by \eqref{eq:twothirds}, or there are at most two points, each with $n=1$.
\end{itemize}
\end{proof}

\begin{proof}[Proof of Theorem~\ref{thm:nu15}]
Combine Propositions~\ref{prop:notcube}, \ref{prop:casea}, \ref{prop:caseb} and~\ref{prop:casec}.
\end{proof}

%======================================================================
\section{The enumeration for \texorpdfstring{$\nu\in\{2,3\}$}{nu in \{2,3\}} and the proof of Theorem B}\label{sec:enum}
%======================================================================

\begin{lemma}\label{lem:constraints}
Let $\Lambda$ be a $6$-candidate with a critical realization of defect $\nu\ge2$. Then:
\begin{enumerate}[label=\textup{(E\arabic*)},leftmargin=3em]
\item $X\in\{40,41\}$;
\item every point has multiplicity $\le7$, and some point has multiplicity $\ge4$;
\item $3\le r_H\le7$, i.e.\ $7\le e(H):=14-r_H\le11$, for every line $H$;
\item $W\ge\max(0,3X-120)$.
\end{enumerate}
\end{lemma}
\begin{proof}
(E1) follows from $\nu=43-X$ and Theorem~\ref{thm:nu15}.

(E2) If all multiplicities were $\le3$, each line would contain exactly $7$ triple points, so $X=35$. (E3) is \eqref{eq:rH}.

(E4) Bojanowski's inequality $t_2+\frac34t_3\ge d+\sum_{r\ge5}(\frac{r^2}4-r)t_r$ \cite{Boj03,Pok21} applies, since $t_r=0$ for $r>10=\frac23d$. As in \cite[Lemma~3.5]{AS}, it is equivalent to $3W\ge5X-4c_2+4$, and $c_2=91-X$.
\end{proof}

\subsection{Generation} The simple rank-$3$ matroids satisfying (E1)--(E4) were generated as in \cite[\S4]{AS}:
\begin{enumerate}[label=\textup{(\arabic*)}]
\item first the configurations of \emph{heavy} points (multiplicity $\ge4$), by canonical augmentation \cite{McK98};
\item then their completions by triple points, by orderly generation with the pruning rules of \cite{AS};
\item finally, isomorph rejection with \texttt{nauty} \cite{MP14}, applied to the incidence graph between lines and points of multiplicity $\ge3$, with the lines as a distinguished colour class.
\end{enumerate}
Write $e(H)=\sum_{p\in H}(m_p-2)=14-r_H$. The pruning rules of \cite{AS} are consequences of three facts:
\begin{itemize}[leftmargin=*]
\item the identity $\sum_H\bigl(e(H)-7\bigr)=(3X-105)-W$;
\item the monotonicity of $X$, $W$ and $e(H)$ along the generation: none of them decreases when a point is added or its multiplicity is raised, and triple points do not change $W$;
\item the bounds $7\le e(H)\le11$ and $W\ge\max(0,3X-120)$ of Lemma~\ref{lem:constraints}.
\end{itemize}
They are valid for any values of these bounds, so the programs of \cite{AS} apply unchanged with the parameters listed in the ancillary README.

Heavy configurations with at most six heavy points were completed one at a time, with a limit of $30$ seconds; the others were completed in batches, without a limit. For $18$ configurations ($5$ with $X=41$ and $13$ with $X=40$) the limit was reached. These were completed by a version of the generator that partitions the search tree at depth $4$ into $50$ parts. The partition is exact by construction. The tree is traversed in the same deterministic order in every part, the nodes at depth $4$ are numbered in that order, and part $j$ explores only the subtrees of the nodes whose number is $\equiv j\pmod{50}$. As a test, the partitioned and the original generator were run on configurations with $111{,}916$ and $74{,}868$ outputs, and the outputs coincide. One of the $13$ deferred configurations with $X=40$ has no completion at all; all $50$ parts of it ran to the end, and its output is empty. Which configurations are deferred depends on the machine; the complete lists do not.

The counts below are numbers of distinct canonical forms over all heavy configurations together. Heavy configurations are pairwise non-isomorphic, and the heavy subconfiguration is an isomorphism invariant, so completions of different heavy configurations are never isomorphic. Consistently with this, the canonical forms from the main run and from the $18$ deferred configurations are pairwise distinct. The SHA-256 checksums of the complete lists are in the ancillary files.

\subsection{Realizability} For every class, the realization space modulo $PGL_3$ was computed exactly over $\QQ$, with the construction plans and elimination of \cite[\S5]{AS}. A plan fixes a projective frame and constructs the remaining lines one by one. For all classes except one for each value of $X$, the construction either forces a degeneracy or ends with at most one free parameter. The two exceptions need a two-parameter plan. The possible outcomes, i.e.\ the statuses returned by the program, are:
\begin{itemize}[leftmargin=*]
\item the construction forces a degeneracy: a constructed point or line vanishes identically (status \texttt{empty-degenerate});
\item the incidence conditions have no solutions, or only degenerate ones, i.e.\ solutions for which the lattice changes (status \texttt{empty});
\item a finite set of points, described by irreducible polynomials over $\QQ$, each giving one Galois orbit;
\item a one-parameter family; this occurred only for four classes with $X=40$, and all members of these families are degenerate;
\item a two-parameter plan, treated by a Gröbner basis computation over $\QQ$ in several frames. Each line through a single constructed point $L_a\cap L_b$ is written $L_a+tL_b$; this chart misses only $L_b$, so it is complete. For both classes every frame gives no solutions or only degenerate ones, and some frame gives the Gröbner basis $\{1\}$.
\end{itemize}
As a positive control, the screening recovers the lattice of the $15$-line subsets of $\cA(G(5,1,3))$ with $\nu=2$ (Example~\ref{ex:g513}), which is known to be realizable and critical.

\begin{theorem}[computer-assisted]\label{thm:enum}
\begin{enumerate}[label=\textup{(\roman*)}]
\item $X=41$ ($\nu=2$): there are $11{,}210{,}809$ isomorphism classes, $11{,}206{,}118$ of them from the main run. Exactly two are realizable over $\CC$. Each has a single Galois orbit of realizations, and both are critical:
\begin{itemize}
\item $(t_3,t_4,t_5)=(5,10,1)$, defined over $\QQ(i)$ by $x^2+2x+2$;
\item $(t_3,t_5,t_6)=(15,1,2)$, defined over $\QQ(\zeta_5)$ by $x^4-3x^3+4x^2-2x+1$.
\end{itemize}
\item $X=40$ ($\nu=3$): there are $15{,}054{,}509$ isomorphism classes, $15{,}050{,}791$ of them from the main run. Exactly four are realizable, each with a single Galois orbit, and none is critical:
\begin{itemize}
\item $(t_3,t_4,t_5)=(7,5,3)$ over $\QQ(\alpha)$ with $\alpha^3=\alpha+1$, with $\mdr\ge8$;
\item $(t_3,t_5)=(10,5)$ over $\QQ(\zeta_5)$, defined by $x^4-3x^3+4x^2-2x+1$, with $\mdr\ge8$;
\item $(t_3,t_4)=(7,11)$ over $\QQ(\sqrt{-3})$, defined by $x^2-x+1$, with $\mdr=7$;
\item $(t_3,t_4)=(10,10)$ over $\QQ(\sqrt{-3})$, defined by $x^2+x+1$, with $\mdr=7$.
\end{itemize}
\end{enumerate}
The values of $\mdr$ were computed exactly over the respective number fields, as the least $q$ with $\dim D_0(\cA)_q>0$, for $q\le7$.
\end{theorem}

\begin{table}[ht]
\centering
\caption{Outcome of the realizability screening, by status of the program. The rows below the first one add up to the number of classes.}\label{tab:screen}
\begin{tabular}{lrr}
\toprule
& $X=41$ & $X=40$\\
\midrule
classes & $11{,}210{,}809$ & $15{,}054{,}509$\\
\midrule
no solutions, or only degenerate ones & $11{,}095{,}780$ & $14{,}934{,}401$\\
construction forces a degeneracy & $115{,}026$ & $120{,}099$\\
one-parameter family, all members degenerate & $0$ & $4$\\
two-parameter plan, only degenerate solutions & $1$ & $1$\\
realizable (one Galois orbit each) & $2$ & $4$\\
\midrule
\quad of which critical & $2$ & $0$\\
\bottomrule
\end{tabular}
\end{table}

\begin{example}\label{ex:g513}
The critical class over $\QQ(\zeta_5)$ is the lattice of the $15$-line subsets of the reflection arrangement
\[
\cA(G(5,1,3))=\{xyz(x^5-y^5)(y^5-z^5)(z^5-x^5)=0\}
\]
with $\nu=2$. These subsets are obtained by deleting one coordinate line, say $x=0$, and two of the lines $y=\zeta_5^jz$ through the opposite coordinate point. All $30$ of them are isomorphic.

For these subsets the canonical foliation is defined by the degree-$6$ basis derivation $\theta$ of the free arrangement $\cA(G(5,1,3))$, whose exponents are $(1,6,11)$. Indeed $\theta$ is tangent to the subset, and $\dim D_0(\cA)_6=1$ by Lemma~\ref{lem:sat}. All $18$ lines are $\theta$-invariant and $\deg F_\theta=18$, so Proposition~\ref{prop:residual} gives $F_\theta=\mathrm{const}\cdot Q_{\cA(G(5,1,3))}$. Hence the residual $R$ is the product of the three omitted lines, which are not concurrent. For the class over $\QQ(i)$, an exact computation over $\QQ(i)$ (ancillary file \texttt{residual\_exact.py}) shows that $R=LM_1M_2$ is a product of three non-concurrent lines: $L$ is $\theta$-invariant, defined over $\QQ(i)$ and not in $\cA$, while $M_1,M_2$ are not invariant, are conjugate over $\QQ(i)$, are defined over $\QQ(\sqrt{-5})$, and meet at a triple point of $\cA$. Since $R$ has multiplicity $2$ at that vertex, (L1) excludes it from $Z$, and $Z$ consists of the two points $L\cap M_1$ and $L\cap M_2$. In both cases $R$ is not a cube, in accordance with \S\ref{ssec:notcube}.
\end{example}

\begin{proof}[Proof of Theorem B]
By Theorem~\ref{thm:induction} with $S=6$ it suffices to prove $(C_6)$, the cases $s\le5$ being part of the proof of Theorem A. Let $\Lambda$ be a $6$-candidate all of whose realizations have $\mdr\ge6$, and suppose it has a critical realization, of defect $\nu$.
\begin{itemize}[leftmargin=*]
\item If $\nu\le1$, Proposition~\ref{prop:nu1}, with Terao's conjecture for free arrangements of $15$ lines with exponents $(7,7)$ \cite{AS}, shows that all realizations are critical.
\item Otherwise $\nu\in\{2,3\}$ by Theorem~\ref{thm:nu15}. By Lemma~\ref{lem:constraints}, $\Lambda$ occurs in the enumeration, so it is one of the two critical classes of Theorem~\ref{thm:enum}. Its realizations form a single Galois orbit, so all of them are critical by Remark~\ref{rem:galois}.
\end{itemize}
\end{proof}

%======================================================================
\section{Consequences for any number of lines}\label{sec:anyd}
%======================================================================

The results of \S\ref{sec:reduction} hold for any number of lines. In this section $d$ is arbitrary.

\begin{proposition}\label{prop:anyd}
Let $S\ge0$ and assume $(C_s)$ for all $s\le S$. Let $\cA$ and $\cB$ be arrangements of $d$ lines with $L(\cA)\cong L(\cB)$.
\begin{enumerate}[label=\textup{(\alph*)}]
\item For every $s\le S$ with $s<K$, $\mdr(\cA)\le s$ if and only if $\mdr(\cB)\le s$.
\item If $\min\{k(\cA),k(\cB)\}\le S$, then $k(\cA)=k(\cB)$.
\item If $\cA$ is free with exponents $(d_1,d_2)$, $d_1\le d_2$, and $d_1\le S+1$, then $\cB$ is free.
\end{enumerate}
\end{proposition}
\begin{proof}
(a) Let $s<K$. By Theorem~\ref{thm:reduction}, $\mdr(\cA)\le s$ if and only if $\mathrm{P}(s')$ holds for $\cA$ for some $s'\le s$. The first clauses of $\mathrm{P}(s')$ depend only on the lattice. The last clause is $[\mdr\le s_0]$ for a sub-arrangement of $2s_0+3$ lines, which may be chosen by a rule depending only on the labelled lattice; it realizes an $s_0$-candidate or a pencil, with $s_0\le s'\le S$. By Theorem~\ref{thm:induction}, $[\mdr\le s_0]$ is constant on $\cR(\Lambda)$ for every $s_0$-candidate $\Lambda$ with $s_0\le S$, and it holds for pencils. So $\mathrm{P}(s')$ has the same truth value for $\cA$ and $\cB$.

(b) Let $k=k(\cA)\le S$, say. If $k<K$, then $\mdr(\cA)=k$, and (a) with $s=k$ and, if $k\ge1$, with $s=k-1$ gives $\mdr(\cB)=k$, so $k(\cB)=k$. If $k=K$, then $K\le S$, and (a) for all $s<K$ gives $\mdr(\cB)\ge K$, so $k(\cB)=K$.

(c) Since $d_1+d_2=d-1$, we have $d_1\le K$, so $s:=d_1-1$ satisfies $s<K$ and $s\le S$. As $\mdr(\cA)=d_1>s$, (a) gives $\mdr(\cB)\ge d_1$. The Chern classes of $\TA$ and $\mathcal T\langle\cB\rangle$ agree, and $c_2=d_1d_2$.

If $d_1=d_2$, then $d=2K+1$, $d_1=K$, and $F=\mathcal T\langle\cB\rangle(K)$ has $c_1(F)=0$, $c_2(F)=K^2-2K^2+K^2=0$ and $H^0(F(-1))=D_0(\cB)_{K-1}=0$. By Lemma~\ref{lem:c1c2zero}, $\cB$ is free.

Let $d_1<d_2$ and $r=\mdr(\cB)$. If $2r\ge d-1$, then $\mathcal T\langle\cB\rangle$ is semistable (proof of Proposition~\ref{prop:gst}), and Bogomolov's inequality gives $4d_1d_2=4c_2\ge(d-1)^2=(d_1+d_2)^2$, which is false. Hence $2r<d-1$, and $c_2\ge r(d-1-r)$ by Proposition~\ref{prop:freecrit}. The function $t\mapsto t(d-1-t)$ is strictly increasing for $t<\frac{d-1}2$, and $r\ge d_1$, so $c_2=d_1(d-1-d_1)$ forces $r=d_1$ and $c_2=r(d-1-r)$. By Proposition~\ref{prop:freecrit}, $\cB$ is free.
\end{proof}

\begin{proof}[Proof of Theorem C]
By Theorems A and B, $(C_s)$ holds for all $s\le6$, and for $s\le5$ without computer assistance (proof of Theorem A). Apply Proposition~\ref{prop:anyd} with $S=6$, respectively $S=5$. For the second statement of (a), note that $2\mdr(\cA)\le d-3$ means $\mdr(\cA)<K$, so $k(\cA)=\mdr(\cA)$, and then $k(\cB)=k(\cA)<K$ gives $\mdr(\cB)=\mdr(\cA)$.
\end{proof}

\subsection{Seventeen and eighteen lines}
Let $\cA$ be free with exponents $(d_1,d_2)$, $d_1\le d_2$, and let $H\in\cA$. By a theorem of Abe \cite[Cor.~1.2(1)]{Abe14}, either $r_H\le d_1+1$ or $r_H=d_2+1$. If $r_H\in\{d_1+1,d_2+1\}$, then $\chi(\cA^H,t)=(t-1)(t-r_H+1)$ divides $\chi(\cA,t)=(t-1)(t-d_1)(t-d_2)$. This is a condition on the lattice, and $\cA^H$ is free since it has rank $2$; so by Abe's division theorem \cite{Abe16} every arrangement with lattice $L(\cA)$ is free. Hence, if some $\cB$ with $L(\cB)\cong L(\cA)$ is not free, every line of $\cA$ contains at most $d_1$ intersection points. This step uses no case of Terao's conjecture for fewer lines; it is part of the numerical reduction of \cite{AS}. Since $d_1\le K$ and $d_1\le7$ is excluded by Theorem C(b), for $d=17$ and $d=18$, where $K=8$, only the exponents $(8,8)$ and $(8,9)$ remain.

\begin{proposition}\label{prop:1718}
\begin{enumerate}[label=\textup{(\arabic*)}]
\item Terao's conjecture holds for $17$ lines if and only if $(C_7)$ holds for the $7$-candidates with $\nu=1$.
\item If $(C_7)$ holds for the $7$-candidates with $\nu\in\{2,3\}$, then Terao's conjecture holds for $18$ lines.
\end{enumerate}
\end{proposition}
\begin{proof}
By Theorems A and B, Proposition~\ref{prop:anyd} applies with $S=6$. For $d=17,18$ we have $K=8$, so $[\mdr\le6]$ is determined by the lattice.

(1) Suppose $(C_7)$ fails for a $7$-candidate $\Lambda$ with $\nu=1$: all realizations have $\mdr\ge7$, some $\cA$ is critical, and some $\cB$ is not. Then $\cA$ is not free, by Lemma~\ref{lem:cand}(b), while $\cB$ is free by Proposition~\ref{prop:nu1}(b). So Terao's conjecture fails for $17$ lines.

Conversely, let $\cA$ be free and $\cB$ not free, with $17$ lines and $L(\cA)\cong L(\cB)$. By the reduction above, $\cA$ has exponents $(8,8)$ and at most $8$ intersection points on every line, so $\Lambda=L(\cA)$ is a $7$-candidate, with $\nu=c_2-63=1$; it is not a pencil, since $d_1>0$. Since $\mdr(\cA)=8$, Proposition~\ref{prop:anyd}(a) shows that all realizations of $\Lambda$ have $\mdr\ge7$. If $\mdr(\cB)\ge8$, then $F=\mathcal T\langle\cB\rangle(8)$ has $c_1(F)=0$, $c_2(F)=64-128+64=0$ and $H^0(F(-1))=D_0(\cB)_7=0$, so $\cB$ would be free by Lemma~\ref{lem:c1c2zero}. So $\mdr(\cB)=7$, $\cB$ is critical and $\cA$ is not, and $(C_7)$ fails for $\Lambda$.

(2) Let $\cA$ be free and $\cB$ not free, with $18$ lines and $L(\cA)\cong L(\cB)$. By the reduction above, $\cA$ has exponents $(8,9)$ and at most $8$ intersection points on every line. By Proposition~\ref{prop:anyd}(a), $\mdr(\cB)\ge7$. If $r=\mdr(\cB)\ge9$, then $2r>d-1$, and Bogomolov's inequality for the semistable bundle $\mathcal T\langle\cB\rangle$ would give $4\cdot72\ge17^2$. If $\mdr(\cB)=8$, then $c_2=72=8\cdot9$ and $\cB$ is free by Proposition~\ref{prop:freecrit}. Hence $\mdr(\cB)=7$.

Here $X=\binom{17}2-72=64$. Since $\sum_pm_p(m_p-2)=2X+\sum_p(m_p-2)=3X-W$, counting incidences gives
\[
\sum_Hr_H=\sum_pm_p=\sum_pm_p(m_p-1)-\sum_pm_p(m_p-2)=d(d-1)-(3X-W)=114+W>6\cdot18 ,
\]
so some line $H$ of $\cA$ has $r_H\in\{7,8\}$. Take the corresponding line of $\cB$, and put $\cA'=\cA\setminus H$, $\cB'=\cB\setminus H$. Their lattices are isomorphic. They have $17$ lines and at most $8$ intersection points on every line, so $\Lambda'=L(\cA')$ is a $7$-candidate; it is not a pencil, since a point of multiplicity $m$ gives $m\le8$ intersection points on a line not through it. Deleting $H$ lowers $c_2=b_2-(d-1)$ by $r_H-1$, so $\nu(\Lambda')=72-(r_H-1)-63=10-r_H\in\{2,3\}$.

By Proposition~\ref{prop:ad}(i), $\mdr(\cA')\in\{7,8\}$. If $\mdr(\cA')=7$, then $d-1-r_H=17-r_H\ge9>7$, and Proposition~\ref{prop:ad}(ii) would give $\mdr(\cA)=7$; so $\mdr(\cA')=8$. On the other hand $\mdr(\cB')\le\mdr(\cB)=7$ by Proposition~\ref{prop:ad}(i), and $\mdr(\cB')\ge7$ by Proposition~\ref{prop:anyd}(a) applied to $\cA',\cB'$. So all realizations of $\Lambda'$ have $\mdr\ge7$, $\cB'$ is critical and $\cA'$ is not, and $(C_7)$ fails for $\Lambda'$.
\end{proof}

By Theorem~\ref{thm:induction}, $(C_7)$ would give \GST\ for at most $18$ lines.

%======================================================================
\section{Remarks and questions}
%======================================================================

\begin{remark}
The proof shows more than $(C_6)$: \emph{every} critical $15$-line arrangement has $\nu\le3$, and those with $\nu\ge2$ are realizations of the two lattices of Theorem~\ref{thm:enum}. Together with $\nu\le1$ for $13$ lines and $\nu=0$ for at most $11$ lines, this suggests the following. For $s=6$ the maximal defect that occurs is $2$, while the bound we prove is $3$.
\end{remark}

\begin{question}
Is $\nu\le s-3$ for every critical arrangement of $2s+3$ lines?
\end{question}

We add an observation that we have not verified systematically and that is not used anywhere. The critical non-free arrangements with $s=7$ ($17$ lines, residual a quartic) that we found in exploratory computer experiments, not included in the ancillary files, have $\nu\le4$. In these examples the residual is the product of two invariant lines and a non-invariant conic, and adding the two lines produces a free arrangement with exponents $(7,11)$. The example of $\cA(G(5,1,3))$ behaves in the same way: deleting the line $x=0$ gives a free arrangement with exponents $(6,10)$, since $x=0$ contains $7$ intersection points, and the critical $15$-line subsets are obtained from it by deleting two $\theta$-invariant lines.

\begin{question}
Is every critical non-free arrangement obtained from a free arrangement with exponents $(s,s+4)$ by deleting two lines that are invariant for the canonical foliation?
\end{question}

For $13$ lines there are examples of this kind. The reflection arrangements $\cA(H_3)$, over $\QQ(\sqrt5)$, and $\cA(G(4,1,3))$, over $\QQ(i)$, have $15$ lines and are free with exponents $(5,9)$. Deleting two lines that meet in a double point gives, in both cases, a critical $5$-candidate with $\nu=1$, whose canonical foliation is defined by the degree-$5$ basis derivation of the reflection arrangement; in particular the two deleted lines are invariant. This was checked in exact arithmetic for every double point (ancillary file \texttt{thirteen.py}). We do not know whether every critical non-free arrangement of $13$ lines arises in this way.

\begin{remark}[on $(C_7)$]
By Proposition~\ref{prop:1718}, the stratum $\nu=1$ of $(C_7)$ is Terao's conjecture for $17$ lines, and the strata $\nu\in\{2,3\}$ imply Terao's conjecture for $18$ lines. Among the strata that are not automatic, it is also the only one in which positive-dimensional realization spaces, and hence jumps of $\mdr$ inside a family, are expected. Indeed, for an $s$-candidate the identity $\sum_H\bigl(e(H)-(s+1)\bigr)=3X-W-d(s+1)\ge0$ gives
\[
2d-8-\sum_p(m_p-2)\ \le\ s-3-2\nu
\]
for the expected dimension of the realization space. The critical locus, where $h^0(\TA(s))>0$, has expected codimension $\nu$, since $\chi(\TA(s))=1-\nu$ and $h^2(\TA(s))=h^0(\TA(s-1))=0$ when $\mdr\ge s$. For $s=7$ the bound is $\le0$ when $\nu\ge2$, but it is $2$ when $\nu=1$. These counts are heuristic: they are not upper bounds for the dimension of the realization space, since the incidence conditions need not be independent, and they do not exclude positive-dimensional realization spaces.

For $s=7$ the residual is a quartic, and for $\nu\ge2$ the method of \S\ref{sec:fifteen} needs a defect bound. New local situations appear: a double conic $R=C^2$ with $C$ a smooth conic, the case $R=L^4$ with jets of higher order, and points of multiplicity $4$ of $R$, at which $Z$ need not be curvilinear; for instance $I_{Z,p}=(x^2,y^2)$ is compatible with $R=x^4-y^4\in I_{Z,p}^2$. The three mechanisms of \S\ref{sec:fifteen} (Riccati, Ziegler, normal budget) work along a line, so they apply to multiple linear components of the residual; a double conic needs a new argument, and the local computations grow with $s$.

The enumeration of \S\ref{sec:enum} does not extend directly either. With the constraints of Lemma~\ref{lem:constraints} adapted to $s=7$ (multiplicities $\le8$, $8\le e(H)\le13$, $X=55$, $W\ge7$), the first stage, \texttt{heavyaug 17 8 8 13 55 7}, produced more than $4\cdot10^7$ heavy configurations before we stopped it, against $96\,828$ for $s=6$ and $X=41$. This computation is not part of any proof in this paper.
\end{remark}

%======================================================================
\appendix
\section{The computer-verified local statements}\label{app:local}
%======================================================================

\subsection*{Setting and exactness}
All local statements in \S\ref{sec:fifteen} concern a point $p=0$ on $\ell=\{x=0\}$ and the affine inflection polynomial $D$ of \eqref{eq:affine}, computed from jets of the local frame coefficients. The jet coefficients are independent indeterminates over $\QQ[\mu,\text{slopes}]$, and we use truncated power series in $x,y$ with coefficients in multivariate polynomial rings over $\QQ$ (\texttt{python-flint} \cite{FLINT}).

\begin{lemma}[exactness]\label{lem:exact}
Suppose the homogeneous parts of the components of $X$ are known up to degree $M$ inclusive. If the components lie in $\mm^2$, then every coefficient of $D$ of total degree $\le M+3$ is computed correctly. If they only lie in $\mm$, the same holds up to degree $M+1$.
\end{lemma}
\begin{proof}
$D$ is a sum of products of three components or their derivatives, with one derivative in total. A monomial of degree $e$ in $D$ therefore involves component monomials of degrees $k_1+k_2+k_3=e+1$ with all $k_i\ge2$ (resp.\ $\ge1$). So every $k_i\le e-3$ (resp.\ $\le e-1$), and this is $\le M$ when $e\le M+3$ (resp.\ $e\le M+1$).
\end{proof}

Only coefficients within these bounds are used, and the jets are chosen large enough for every order that occurs; \texttt{verify\_all.py} asserts the bound of Lemma~\ref{lem:exact} for every coefficient it reads. Each induction step below has the following form: assuming that the lower pure coefficients vanish, the displayed coefficient of $D$ is computed and found to be \emph{exactly} a nonzero factor times the next pure coefficient. The factor is a nonzero rational constant or a nonzero multiple of the nondegeneracy factor. The script \texttt{verify\_all.py} performs all $87$ steps and asserts each of them. It also asserts the displayed low-order coefficients in (A2) and (C2) and the branch equations in (C3) for $m=3,4$. Its log is \texttt{verify\_all.txt}. In the vertex cases the slopes beyond the normalized $0,1$ are kept as indeterminates $s_3,s_4,\dots$, and no step depends on them.

\begin{enumerate}[label=(\Alph*0),leftmargin=3em]
\item[(A0)] (L3), branch $g_1=\mu f_1$: if $n_p\ge3$, then $\inn_3(D/u)$ is $-\mu v^3(\dots)$ or $-\mu u^3(\dots)$, according to the tangent of $Z_p$ (\texttt{typeS\_A.py}).
\item[(A1)] The identities of Lemma~\ref{lem:riccati}, with $f_i,g_i$ arbitrary functions of $v$, and the vertex formula of Lemma~\ref{lem:vertexformula} for $m=2,3$ with symbolic slope.
\item[(A2)] Case (a), type S: for general linear parts, $[x^2y]=-b_{01}\Delta$, $[x^3]=-b_{10}\Delta$, and $[x^4]$ is a combination of $b_{10},b_{01}$. At $p_0$, with $c_0=ya_0$: $[x^3y^{j-1}]=-j\,a_{0,j}$ for $2\le j\le7$.
\item[(A3)] Case (a), vertex $q\neq p_0$: coefficient of $x^3y^{j+m-1}$ is $-m(m+j-2)\,c_{0,j}$, for $4\le m\le7$, $2\le j\le m-2$. The script works in a chart in which $p_0=(0,-\delta)$ is finite and uses the Riccati relation $g_0=(y+\delta)f_0$ in the form $\delta a_0=-b_0y^{m-2}\bigl(my+(m-1)\delta\bigr)$, with $1/\delta$ as an indeterminate. The factors found do not involve $1/\delta$, so the case $1/\delta=0$, i.e.\ $p_0$ at infinity as in the proof of Lemma~\ref{lem:casea}(b), is included.
\item[(A4)] Case (a), vertex at $p_0$, $\mu$-branch: coefficient of $x^3y^{j+m-1}$ is $-m(2m-1-j)\mu\,a_{0,j}$, for $3\le m\le7$, $2\le j\le4$.
\item[(A5)] Case (a), double point at $p_0$, B-branch: coefficient of $x^3y^{j+1}$ is $-2(j+1)a_{0,j}$, $2\le j\le7$.
\item[(C1)] Case (c), type O, saddle-node case: the coefficient expansions used in Lemma~\ref{lem:casecO}; the coefficient of $x^2y^j$ is $3a_{0,j}$ for $2\le j\le7$. The nilpotent case of Lemma~\ref{lem:casecO} is the computation in its proof, not part of this item.
\item[(C2)] Case (c), type S: the low-order coefficients
\begin{align*}
[xy^2]&=a_{10}(a_{01}b_{10}-a_{10}b_{01}),\qquad [x^2y^2]=3a_{10}(a_{01}b_{20}-a_{20}b_{01}+b_{01}b_{10}),\\
[x^3y]&=a_{10}(a_{10}b_{20}-a_{20}b_{10}+b_{10}^2),
\end{align*}
and then, for $2\le j\le7$, $[xy^{j+1}]=-c_{0,j}$ and $[x^2y^{j+1}]=3(c_{20}+\mu^2)a_{0,j}$.
\item[(C3)] Case (c), vertices.
\begin{itemize}
\item $m=2$, $\mu$-branch: $[xy^{j+3}]=4\mu\,c_{0,j}$ and $[x^2y^{j+2}]=6\mu(\mu^2-1)a_{0,j}$.
\item $m=2$, B-branch: $-4a_{0,j}$ and $6c_{0,j}$; here $2\le j\le7$.
\item $m=3$: the branch equations stated in Lemma~\ref{lem:casec}(c); then, with $c_{0,j}$ eliminated, $[xy^{j+5}]=-\frac{27}2\mu^2a_{0,3}$, $-\frac{54}5\mu^2a_{0,4}$, $-9\mu^2a_{0,5}$.
\item $m=4$: the relations stated in Lemma~\ref{lem:casec}(d).
\end{itemize}
\item[(C4)] $R=L_1L_2^2$ with $L_2$ not invariant, a point of type S in $Z$: the coefficient of $x^2y$ in $D$ is identically $0$.
\item[(B1)] Multiarrangements (an independent check of Lemma~\ref{lem:multiarr}): for $r=3$ the seven multiplicity types give $\dim D_6=0$. For $r=4$ the gcd of the maximal minors for the $16$ types is one of $\lambda^a(\lambda-1)^b$, with $a,b\in\{0,1,2\}$, or a nonzero constant.
\end{enumerate}

\subsection*{Enumeration data}
Heavy configurations: \texttt{heavyaug}; completion: \texttt{lsgen3} and \texttt{lsgen3s} (partitioned); screening: \texttt{screen\_s6.py}, using \texttt{exact.py} of \cite{AS}; values of $\mdr$: \texttt{mdr\_s6.py}; two-parameter plans: \texttt{cost2.py}. The ancillary files contain the programs, the screening logs, the realizations of the six realizable classes, and SHA-256 checksums of the complete enumeration output, which is available on request (about $2.9$\,GB) and can be regenerated deterministically (see the README).

\subsection*{Declarations}
\textbf{Competing interests.} The authors have no financial interests to disclose.

\textbf{Data availability.} The code and the results of the computations are available as ancillary files of the arXiv version of this paper (Sections~\ref{sec:enum} and~\ref{app:local}).

\textbf{Use of generative AI.} Claude (Anthropic) was used for literature search, programming, and drafting the text; ChatGPT (OpenAI) for auditing and editorial assistance. After using these tools, the authors reviewed and edited the content and take full responsibility for the content of the publication.


\begin{thebibliography}{CHMN18}

\bibitem[Abe13]{Abe13}
T.~Abe, \emph{Chambers of $2$-affine arrangements and freeness of $3$-arrangements}, J.\ Algebraic Combin.\ \textbf{38} (2013), 65--78; arXiv:1005.5276 (as \emph{Exponents of $2$-multiarrangements and freeness of $3$-arrangements}).

\bibitem[Abe14]{Abe14}
T.~Abe, \emph{Roots of characteristic polynomials and intersection points of line arrangements}, J.\ Singul.\ \textbf{8} (2014), 100--116.

\bibitem[Abe16]{Abe16}
T.~Abe, \emph{Divisionally free arrangements of hyperplanes}, Invent.\ Math.\ \textbf{204} (2016), 317--346.

\bibitem[AD18]{AD18}
T.~Abe, A.~Dimca, \emph{Splitting types of bundles of logarithmic vector fields along plane curves}, Internat.\ J.\ Math.\ \textbf{29} (2018), no.~8, 1850055.

\bibitem[ADS21]{ADS21}
T.~Abe, A.~Dimca, G.~Sticlaru, \emph{Addition--deletion results for the minimal degree of logarithmic derivations of hyperplane arrangements and maximal Tjurina line arrangements}, J.\ Algebraic Combin.\ \textbf{54} (2021), 739--766.

\bibitem[AS]{AS}
C.~Anghel, R.~Stan, \emph{Terao's freeness conjecture for arrangements of $15$ and $16$ lines}, arXiv:2610.25868, 2026.

\bibitem[BBJKL21]{BBJKL21}
M.~Barakat, R.~Behrends, C.~Jefferson, L.~Kühne, M.~Leuner, \emph{On the generation of rank $3$ simple matroids with an application to Terao's freeness conjecture}, SIAM J.\ Discrete Math.\ \textbf{35} (2021), 1201--1223.

\bibitem[BK23]{BK23}
M.~Barakat, L.~Kühne, \emph{Computing the nonfree locus of the moduli space of arrangements and Terao's freeness conjecture}, Math.\ Comp.\ \textbf{92} (2023), 1431--1452.

\bibitem[Boj03]{Boj03}
R.~Bojanowski, \emph{Zastosowania uogólnionej nierówności Bogomolova--Miyaoka--Yau}, Master's thesis, University of Warsaw, 2003.

\bibitem[CHMN18]{CHMN18}
D.~Cook II, B.~Harbourne, J.~Migliore, U.~Nagel, \emph{Line arrangements and configurations of points with an unexpected geometric property}, Compos.\ Math.\ \textbf{154} (2018), 2150--2194.

\bibitem[Dim24]{Dim24}
A.~Dimca, \emph{On free curves and related open problems}, Rev.\ Roumaine Math.\ Pures Appl.\ \textbf{69} (2024), no.~2, 129--150; arXiv:2312.07591.

\bibitem[DIM19]{DIM19}
A.~Dimca, D.~Ibadula, A.~Măcinic, \emph{Freeness for $13$ lines arrangements is combinatorial}, Discrete Math.\ \textbf{342} (2019), 2445--2453.

\bibitem[DP26a]{DP26}
A.~Dimca, P.~Pokora, \emph{A nine-line counterexample to a conjecture on the minimal degree of Jacobian relations}, arXiv:2607.01985, 2026.

\bibitem[DP26b]{DPgen}
A.~Dimca, P.~Pokora, \emph{A generalized Terao Conjecture for line arrangements}, arXiv:2609.33584v2, 2026.

\bibitem[EGH96]{EGH96}
D.~Eisenbud, M.~Green, J.~Harris, \emph{Cayley--Bacharach theorems and conjectures}, Bull.\ Amer.\ Math.\ Soc.\ (N.S.) \textbf{33} (1996), 295--324.

\bibitem[FLINT]{FLINT}
The FLINT team, \emph{FLINT: Fast Library for Number Theory}, via python-flint 0.9.0, \url{https://flintlib.org}.

\bibitem[FV14]{FV14}
D.~Faenzi, J.~Vallès, \emph{Logarithmic bundles and line arrangements, an approach via the standard construction}, J.\ London Math.\ Soc.\ (2) \textbf{90} (2014), 675--694.

\bibitem[McK98]{McK98}
B.~D.~McKay, \emph{Isomorph-free exhaustive generation}, J.\ Algorithms \textbf{26} (1998), 306--324.

\bibitem[MP14]{MP14}
B.~D.~McKay, A.~Piperno, \emph{Practical graph isomorphism, II}, J.\ Symbolic Comput.\ \textbf{60} (2014), 94--112.

\bibitem[OSS]{OSS}
C.~Okonek, M.~Schneider, H.~Spindler, \emph{Vector Bundles on Complex Projective Spaces}, Progr.\ Math.\ \textbf{3}, Birkhäuser, Boston, 1980.

\bibitem[OT92]{OT92}
P.~Orlik, H.~Terao, \emph{Arrangements of Hyperplanes}, Grundlehren Math.\ Wiss.\ \textbf{300}, Springer, Berlin, 1992.

\bibitem[Per01]{Per01}
J.~V.~Pereira, \emph{Vector fields, invariant varieties and linear systems}, Ann.\ Inst.\ Fourier \textbf{51} (2001), 1385--1405.

\bibitem[Pok21]{Pok21}
P.~Pokora, \emph{Hirzebruch-type inequalities viewed as tools in combinatorics}, Electron.\ J.\ Combin.\ \textbf{28} (2021), \#P1.9.

\bibitem[Sai80]{Sai80}
K.~Saito, \emph{Theory of logarithmic differential forms and logarithmic vector fields}, J.\ Fac.\ Sci.\ Univ.\ Tokyo Sect.\ IA Math.\ \textbf{27} (1980), 265--291.

\bibitem[Wak07]{Wak07}
A.~Wakamiko, \emph{On the exponents of $2$-multiarrangements}, Tokyo J.\ Math.\ \textbf{30} (2007), 99--116.

\bibitem[Zie89]{Zie89}
G.~M.~Ziegler, \emph{Multiarrangements of hyperplanes and their freeness}, in: Singularities (Iowa City, IA, 1986), Contemp.\ Math.\ \textbf{90}, Amer.\ Math.\ Soc., Providence, RI, 1989, 345--359.

\bibitem[Zie89a]{Zie89a}
G.~M.~Ziegler, \emph{Combinatorial construction of logarithmic differential forms}, Adv.\ Math.\ \textbf{76} (1989), 116--154.

\end{thebibliography}
\end{document}